\documentclass{article}
\usepackage[dvips]{graphicx}
\usepackage{amsmath}
\usepackage{amsfonts}
\usepackage{amssymb}
\newtheorem{theorem}{Theorem}[subsection]

\newtheorem{definition}[theorem]{Definition}
\newtheorem{example}[theorem]{Example}

\newtheorem{lemma}[theorem]{Lemma}

\newtheorem{proposition}[theorem]{Proposition}
\newtheorem{remark}[theorem]{Remark}

\newenvironment{proof}[1][Proof]{\textbf{#1.} }{\ \rule{0.5em}{0.5em}}

\begin{document}

\title{An algorithm for computing the global basis of a finite dimensional
irreducible $U_{q}(so_{2n+1})$ or $U_{q}(so_{2n})$-module}
\author{C\'{e}dric Lecouvey\\lecouvey@math.unicaen.fr}
\date{}
\maketitle
\begin{abstract}
We describe a simple algorithm for computing the canonical basis of any
irreducible finite-dimensional $U_{q}(so_{2n+1})$ or $U_{q}(so_{2n})$-module.
\end{abstract}

\section{ Introduction}

The quantum algebra $U_{q}(\frak{g})$ associated to a semisimple Lie algebra
$\frak{g}$ is the $q$-analogue introduced by Drinfeld and Jimbo of its
universal enveloping algebra $U(\frak{g})$. Kashiwara \cite{Ka} and Lusztig
\cite{Lut} have discovered a distinguished basis of $U_{q}^{-}(\frak{g})$
which projects onto a global crystal basis (Kashiwara) or canonical basis
(Lusztig) of each simple finite-dimensional $U_{q}(\frak{g})$-module. When $q$
tends to $1$ this basis yields in particular a canonical basis of the
corresponding $U(\frak{g})$-module.

In this article, we consider the orthogonal case $g=so_{2n+1}$ and $g=so_{2n}%
$. For each dominant weight $\lambda$ there exists a unique irreducible
finite-dimensional $U_{q}(g)$-module $V(\lambda)$ of highest weight $\lambda$.
The aim of this article is to describe a simple algorithm for computing the
global crystal basis of $V(\lambda)$ analogous to those given in \cite{L-T}
and \cite{L1} respectively for the irreducible $U_{q}(sl_{n})$ and
$U_{q}(sp_{2n})$-modules.

For type $B,$ we denote by $\{\Lambda_{1}^{B},...,\Lambda_{n}^{B}\}$ the set
of fundamentals weights of $U_{q}(so_{2n+1}).$ Set $\omega_{p}^{B}=\Lambda
_{p}^{B}$ for $p=1,...,n-1$ and $\omega_{n}^{B}=2\Lambda_{n}^{B}.$ Then
$V^{B}(\lambda)$ is isomorphic to an irreducible component of a tensor power
of the vector representation if and only if $\lambda$ belongs to the subset
$\Omega_{+}^{B}$ of dominant weights which can be decomposed on the basis
$\{\omega_{p}^{B},$ $p=1,...,n\}$.\ In particular the spin representation
$V^{B}(\Lambda_{n})$ is not isomorphic to an irreducible component of a tensor
power of the vector representation.\ Every dominant weight of $U_{q}%
(so_{2n+1})$ may be written $\lambda=\Lambda+\lambda^{\prime}$ with
$\Lambda\in\{0,\Lambda_{n}^{B}\}$ and $\lambda^{\prime}\in\Omega_{+}^{B}.$
Then our method is as follows.

\noindent For any $\omega_{p}^{B},$ we realize first $V(\omega_{p}^{B})$ as a
subrepresentation of a $U_{q}(so_{2n+1})$-module $W(\omega_{p}^{B})$ whose
basis $\{v_{C}\}$ has a natural indexation in terms of column shaped Young
tableaux. This representation $W(\omega_{p}^{B})$ may be regarded as a
$q$-analogue of the $p$-th exterior power of the vector representation of
$U_{q}(so_{2n+1})$. It has been inspired by a similar construction for
quantized affine Lie algebras obtained by Kashiwara, Miwa, Petersen and Yung
in \cite{KMPY}. More precisely we have%
\[
W(\omega_{p}^{B})=V(\Lambda_{1}^{B})^{\otimes p}/N_{p}^{B}%
\]
where $N_{p}^{B}=\underset{i=0}{\overset{p-2}{\sum}}V(\Lambda_{1}%
^{B})^{\otimes i}\otimes N^{B}\otimes V(\Lambda_{1}^{B})^{\otimes(p-2-i)}$ and
$N^{B}$ is the submodule of $V(\Lambda_{1}^{B})^{\otimes2}$ isomorphic to
$V(2\Lambda_{1}^{B}).\;$This construction may be generalized to types $C$ and
$D.\;$For type $C,$ the $U_{q}(sp_{2n})$-module $W(\Lambda_{p}^{C})$ is
precisely that introduced in \cite{L1} by describing the action of Chevalley's
operators on a suitable basis.\ By using a general algorithm due to Marsh
\cite{Ma}, we can compute the decomposition of the canonical basis of
$V(\omega_{p}^{B})$ on the basis $\{v_{C}\}.\;$Note that a $q$-analogue
$\mathcal{W(}\omega_{p}^{B})$ of the $p$-th exterior power of the vector
representation already exists.\ It has been introduced by Jing, Misra and
Okado in \cite{JMO}.\ The module $\mathcal{W(}\omega_{p}^{B})$ is defined by
turning $N$ into $\mathcal{N=}N%
%TCIMACRO{\tbigoplus }%
%BeginExpansion
{\textstyle\bigoplus}
%EndExpansion
V(0)$ into the definition of $W(\omega_{p}^{B}).\;\mathcal{W(}\omega_{p}^{B})$
is irreducible but seems not to allow a generalization of the algorithm
described in \cite{L-T} and \cite{L1}.

\noindent Next we embed $V^{B}(\lambda)$ in
\begin{gather*}
W^{B}(\lambda)=W(\omega_{1}^{B})^{\otimes\mu_{1}}\otimes\cdot\cdot\cdot\otimes
W(\omega_{n}^{B})^{\otimes\mu_{n}\text{ }}\text{if }\lambda=\overset
{n}{\underset{i=1}{\sum}}\mu_{i}\omega_{i}^{B}\in\Omega_{+}^{B},\\
W^{B}(\lambda)=V(\Lambda)\otimes W^{B}(\lambda^{\prime})\text{ otherwise.}%
\end{gather*}
The tensor product of the crystal bases of the modules $W(\omega_{i}^{B})$
occurring in $W^{B}(\lambda)$ is a natural basis $\{v_{\tau}\}$ of
$W^{B}(\lambda)$ indexed by combinatorial objects $\tau$ called tabloids. Then
we obtain an intermediate basis of $V^{B}(\lambda)$ fixed by the involution
$q\longmapsto q^{-1}$ and such that the transition matrix from this basis to
the global crystal basis of $V^{B}(\lambda)$ is unitriangular. Finally we
compute the expansion of the canonical basis on the basis $\{v_{\tau}\}$ via
an elementary algorithm. Note that the coefficients of this expansion are
integral that is belong to $\mathbb{Z}[q,q^{-1}].$

For type $D$ the method is essentially the same up to minor modifications due
to the existence of the two spin representations.

\bigskip

\noindent\textbf{Notation \ }In the sequel, we often write $B$ and $D$ instead
of $B_{n}$ and $D_{n}$ to simplify the notation. Moreover, we frequently
define similar objects for types $B$ and $D$. When they are related to type
$B$ (respectively $D$), we attach to them the label $^{B}$ (respectively the
label $^{D}$). To avoid cumbersome repetitions, we sometimes omit the labels
$^{B}$ and $^{D}$ when our statements are true for the two types.

\section{Background}

In this section we briefly review the basic facts that we shall need
concerning the representation theory of $U_{q}(\frak{so}_{2n+1})$ and
$U_{q}(\frak{so}_{2n})$ and the notions of crystal basis and canonical basis
of their representations. The reader is referred to \cite{Ch-Pr}, \cite{Kang},
\cite{Jan}, \cite{Ka1} and \cite{Ka2} for more details.

\subsection{The quantum enveloping algebras $U_{q}(\frak{so}_{2n+1})$ and
$U_{q}(\frak{so}_{2n})$}

Given a fixed indeterminate $q$ set%
\begin{gather*}
q_{i}^{B}=\left\{
\begin{tabular}
[c]{l}%
$q^{2}$ if $i\neq n$\\
$q$ if $i=n$%
\end{tabular}
\right.  \text{, and }q_{i}^{D}=q\text{ for any }i\in\{1,...,n\},\\
\lbrack m]_{i}=\frac{q_{i}^{m}-q_{i}^{-m}}{q_{i}-q_{i}^{-1}}\text{ and
}[m]_{i}!=[m]_{i}[m-1]_{i}\cdot\cdot\cdot\lbrack1]_{i}.
\end{gather*}
If $\frak{g}$ is one of the two algebras $\frak{so}_{2n+1}$ or $\frak{so}%
_{2n},$ the quantized enveloping algebras $U_{q}(\frak{g})$ is the associative
algebra over $\mathbb{Q}(q)$ generated by $e_{i},f_{i},t_{i},t_{i}^{-1},$
$i=1,...,n$, subject to relations which depend on the coefficients of the
Cartan matrix of $\frak{g}.$.\ Note that for any $i\in I,$ the subalgebra of
$U_{q}(\frak{g})$ generated by $e_{i},f_{i}$ and $t_{i}$ is isomorphic to
$U_{q}(\frak{sl}_{2})$ the quantum enveloping algebra associated to
$\frak{sl}_{2}.$

The representation theory of $U_{q}(\frak{g})$ is closely parallel to that of
$\frak{g}$. The weight lattice $P$ of $U_{q}(\frak{g})$ is the $\mathbb{Z}%
$-lattice generated by the fundamentals weights $\Lambda_{1},...,\Lambda_{n}.$
Write $P_{+}$ for the set of dominant weights of $U_{q}(\frak{g})$.\ We denote
by $V(\lambda)$ the irreducible finite dimensional $U_{q}(\frak{g})$-module
with highest weight $\lambda\in P^{+}.$

Given two $U_{q}(\frak{g})$-modules $M$ and $N$, we can define a structure of
$U_{q}(\frak{g})$-module on $M\otimes N$ by putting:%
\begin{gather}
t_{i}(u\otimes v)=t_{i}u\otimes t_{i}v,\label{tensor1}\\
e_{i}(u\otimes v)=e_{i}u\otimes t_{i}^{-1}v+u\otimes e_{i}v,\label{tensor2}\\
f_{i}(u\otimes v)=f_{i}u\otimes v+t_{i}u\otimes f_{i}v. \label{tensor3}%
\end{gather}

In the sequel we need the following general lemma (see \cite{Jan} p.32).\ Let
$V(l)$ be the irreducible $U_{q}(sl_{2})$-module of dimension $l+1$.

\begin{lemma}
\label{lem_sl2}Consider $v_{r}\in V(r)$ and $v_{s}\in V(s)$.\ Set
$t(v_{r})=q^{a}v_{r}$ with $a\in\mathbb{Z}$. Then for any integer $r$ one has:%
\[
f^{(m)}(v_{r}\otimes v_{s})=\overset{m}{\underset{k=0}{\sum}}q^{(m-k)(a-k)}%
\ f^{(k)}(v_{r})\otimes f^{(m-k)}(v_{s}).
\]
\end{lemma}

\subsection{Crystal basis and crystal graph of $U_{q}(\frak{g})$-modules}

The reader is referred to \cite{Jan} and \cite{Ka2} for basic definitions on
crystal bases and crystal graphs. Given $(L,B)$ and $(L^{\prime},B^{\prime})$
two crystal bases of the finite-dimensional $U_{q}(\frak{g})$-modules $M$ and
$M^{\prime}$, $(L\otimes L^{\prime},$ $B\otimes B^{\prime})$ with $B\otimes
B^{\prime}=\{b\otimes b^{\prime};$ $b\in B,b^{\prime}\in B^{\prime}\}$ is a
crystal basis of $M\otimes M^{\prime}$. The action of $\widetilde{e}_{i}$ and
$\widetilde{f}_{i}$ on $B\otimes B^{\prime}$ is given by:%

\begin{align}
\widetilde{f_{i}}(u\otimes v)  &  =\left\{
\begin{tabular}
[c]{c}%
$\widetilde{f}_{i}(u)\otimes v$ if $\varphi_{i}(u)>\varepsilon_{i}(v)$\\
$u\otimes\widetilde{f}_{i}(v)$ if $\varphi_{i}(u)\leq\varepsilon_{i}(v)$%
\end{tabular}
\right. \label{TENS1}\\
&  \text{and}\nonumber\\
\widetilde{e_{i}}(u\otimes v)  &  =\left\{
\begin{tabular}
[c]{c}%
$u\otimes\widetilde{e_{i}}(v)$ if $\varphi_{i}(u)<\varepsilon_{i}(v)$\\
$\widetilde{e_{i}}(u)\otimes v$ if$\varphi_{i}(u)\geq\varepsilon_{i}(v)$%
\end{tabular}
\right.  \label{TENS2}%
\end{align}
where $\varepsilon_{i}(u)=\max\{k;\widetilde{e}_{i}^{k}(u)\neq0\}$ and
$\varphi_{i}(u)=\max\{k;\widetilde{f}_{i}^{k}(u)\neq0\}$. For any $\lambda\in
P_{+},$ write $B(\lambda)$ for the crystal graph of $V(\lambda).$

\subsection{The spin representations}

Now we review Kashiwara-Nakashima's description \cite{KN} of the spin
representation $V(\Lambda_{n}^{B}),$ $V(\Lambda_{n}^{D})$ and $V(\Lambda
_{n-1}^{D})$ of $U_{q}(\frak{so}_{2n+1})$ and $U_{q}(\frak{so}_{2n}).$ It is
based on the notion of spin columns. To avoid confusion between these new
columns and the classical columns of a tableau that we will define in
\ref{sub_sec_combi}, we will follow Kashiwara-Nakashima's convention and
consider spin columns as column shape diagrams of width $1/2$. Such diagrams
will be called spin diagrams.

\begin{definition}
A spin column $\frak{C}$ of height $n$ is a filling of a spin diagram by $n$
letters $\{x_{1}\prec\cdot\cdot\cdot\prec x_{n}\}$ of the totally ordered set
$\{1\prec\cdot\cdot\cdot\prec n\prec\overline{n}\prec\cdot\cdot\cdot
\prec\overline{1}\}$ such that the letters are increasing from top to
bottom.\ The set of spin columns of length $n$ will be denoted $SP_{n}$.
\end{definition}

\noindent To describe the two spin representations of $U_{q}(\frak{so}_{2n})$
we split $SP_{n}$ in two parts. Set
\begin{align*}
SP_{n}^{+}  &  =\left\{  \frak{C}\in SP_{n};\text{ the number of barred
letters in }\frak{C}\text{ is even}\right\}  ,\\
SP_{n}^{-}  &  =\left\{  \frak{C}\in SP_{n};\text{ the number of barred
letters in }\frak{C}\text{ is odd}\right\}  .
\end{align*}
In \cite{KN}, Kashiwara and Nakashima consider the module $V(\Lambda_{n}^{B})$
as a $2^{n}$ dimension space vector of basis $\frak{B=}\{v_{\frak{C}},$
$\frak{C}\in SP_{n}\}$.\ They give the action of Chevalley's operators on each
vector $v_{\frak{C}}$.\ For this action the vector $v_{\frak{C}_{n}}$ (with
$\frak{C}_{n}$ the spin column containing only unbarred letters) is of highest
weight. In the sequel we only use the action of the operators $f_{i}$'s:%
\begin{align*}
f_{i}(v_{\frak{C}})  &  =\left\{
\begin{tabular}
[c]{l}%
$v_{\frak{C}^{\prime}}$ with $\frak{C}^{\prime}=\frak{C}-\{i,\overline
{i+1}\}+\{i+1,\overline{i}\}$ if $i,\overline{i+1}\in\frak{C}$\\
$0$ otherwise
\end{tabular}
\right.  \text{ for }i=1,...,n-1,\\
f_{n}(v_{\frak{C}})  &  =\left\{
\begin{tabular}
[c]{l}%
$v_{\frak{C}^{\prime}}$ with $\frak{C}^{\prime}=\frak{C}-\{n\}+\{\overline
{n}\}$ if $n\in\frak{C}$\\
$0$ otherwise
\end{tabular}
\right.  .
\end{align*}
Similarly the modules $V(\Lambda_{n}^{D})$ and $V(\Lambda_{n-1}^{D})$ can be
respectively regarded as $2^{n-1}$ dimension spaces vector of basis
$\frak{B}_{+}\frak{=}\{v_{\frak{C}},$ $\frak{C}\in SP_{n}^{+}\}$ and
$\frak{B}_{-}\frak{=}\{v_{\frak{C}},$ $\frak{C}\in SP_{n}^{-}\}.$ For the
action defined in \cite{KN}, $v_{\frak{C}_{n}}$ and $v_{\frak{C}_{n-1}}$ (with
$\frak{C}_{n-1}$ the spin column containing only $\overline{n}$ as barred
letter) are respectively highest weight vectors of $V(\Lambda_{n}^{D})$ and
$V(\Lambda_{n-1}^{D}).\;$The action of the operators $f_{i},$ $i=1,...,n-1$ is
the same as in $V(\Lambda_{n}^{B})$ and we have:%
\[
f_{n}(v_{\frak{C}})=\left\{
\begin{tabular}
[c]{l}%
$v_{\frak{C}^{\prime}}$ with $\frak{C}^{\prime}=\frak{C}-\{n,n-1\}+\{\overline
{n-1},\overline{n}\}$ if $n,n-1\in\frak{C}$\\
$0$ otherwise
\end{tabular}
\right.  .
\]
Note that the actions of the $e_{i}$'s and $t_{i}$'s on the spin
representations may be given in a similar simple way.

\subsection{Combinatorics of crystal graphs \label{sub_sec_combi}}

In this paragraph we recall Kashiwara-Nakashima's realization of $B(\lambda)$
\cite{KN}. It is based on the notion of orthogonal tableau analogous to Young
tableau for type $A$.

The crystal graphs of the vectors representations of $U_{q}(\frak{so}_{2n+1})$
and $U_{q}(\frak{so}_{2n})$ are respectively:%
\begin{equation}
B(\Lambda_{1}^{B}):1\overset{1}{\rightarrow}2\cdot\cdot\cdot\rightarrow
n-1\overset{n-1}{\rightarrow}n\overset{n}{\rightarrow}0\overset{n}%
{\rightarrow}\overline{n}\overset{n-1}{\rightarrow}\overline{n-1}\overset
{n-2}{\rightarrow}\cdot\cdot\cdot\rightarrow\overline{2}\overset
{1}{\rightarrow}\overline{1} \label{vect_B}%
\end{equation}
and%
\begin{equation}
B(\Lambda_{1}^{D}):1\overset{1}{\rightarrow}2\overset{2}{\rightarrow}%
\cdot\cdot\cdot\overset{n-3}{\rightarrow}n-2\overset{n-2}{\rightarrow}%
\begin{tabular}
[c]{c}%
$\overline{n}$ \ \ \\
\ \ $\overset{n}{\nearrow}$ $\ \ \ \overset{n-1}{\text{ \ }\searrow}$ \ \ \ \\
$n-1\ \ \ \ \ \ \ \ \ \ \overline{n-1}$\\
\ $\underset{n-1}{\searrow}$ \ \ \ $\underset{n}{\nearrow}$ \ \ \ \\
$n$ \
\end{tabular}
\overset{n-2}{\rightarrow}\overline{n-2}\overset{n-3}{\rightarrow}\cdot
\cdot\cdot\overset{2}{\rightarrow}\overline{2}\overset{1}{\rightarrow
}\overline{1}. \label{vect_D}%
\end{equation}
Set $G_{n}^{B}=\underset{l\geq0}{\bigoplus(}B^{B}(\Lambda_{1}))^{\otimes l}$
and $G_{n}^{D}=\underset{l\geq0}{\bigoplus(}B^{D}(\Lambda_{1}))^{\otimes l}%
$.\ For any $b\in G_{n}$, let $B(b)$ be the connected componet of $G_{n}$
containing $b.$

\noindent The vertices $b_{p}=1\otimes2\cdot\cdot\cdot\otimes p\in G_{n}$,
$p=1,...,n,$ and the vertex $b_{\overline{n}}=1\otimes2\cdot\cdot\cdot
\otimes(n-1)\otimes\overline{n}\in G_{n}^{D}$ are highest weight vertices.
Consider $\lambda=\overset{n}{\underset{i=1}{\sum}}\lambda_{i}\Lambda_{i}$ a
dominant weight.\ It is not always possible to realize $B(\lambda)$ as a
sub-crystal of $G_{n}.$ Indeed, in $G_{n}^{B}$ the vertex $b_{n}$ is of
highest weight $\omega_{n}^{D}=2\Lambda_{n}$. Hence $B(\lambda)$ is a
sub-crystal of $G_{n}^{B}$ only if $\lambda_{n}$ is even.\ Then $B(\lambda)$
is identified with $B(b_{1}^{\otimes\lambda_{1}}\otimes\cdot\cdot\cdot\otimes
b_{n}^{\lambda_{n}/2})$.\ We denote by $\Omega_{+}^{B}$ the sub-set of such
dominant weights. Similarly, in $G_{n}^{D},$ the vertices $b_{n}$,
$b_{\overline{n}}$ and $b_{n-1}$ are respectively of highest weight
$\omega_{n}^{D}=2\Lambda_{n},$ $\overline{\omega}_{n}^{D}=2\Lambda_{n-1}$ and
$\omega_{n-1}^{D}=\Lambda_{n}+\Lambda_{n-1}$.\ Hence, $B(\lambda)$ is a
sub-crystal of $G_{n}^{D}$ if and only if $\lambda_{n}$ and $\lambda_{n-1}$
have the same parity (i.e. $\lambda_{n}-\lambda_{n-1}=0$ $\operatorname{mod}%
2$). Then $B(\lambda)$ is identified with $B(b_{1}^{\otimes\lambda_{1}}%
\otimes\cdot\cdot\cdot\otimes b_{n-1}^{(\lambda_{n}-\lambda_{n-1})/2}\otimes
b_{n})$ if $\lambda_{n}>\lambda_{n-1},$ with $B(b_{1}^{\otimes\lambda_{1}%
}\otimes\cdot\cdot\cdot\otimes b_{n-1}^{(\lambda_{n-1}-\lambda_{n})/2}\otimes
b_{\overline{n}})$ if $\lambda_{n}<\lambda_{n-1}$ and with $B(b_{1}%
^{\otimes\lambda_{1}}\otimes\cdot\cdot\cdot\otimes b_{n-1}^{(\lambda
_{n-1}-\lambda_{n})/2})$ otherwise.\ We denote by $\Omega_{+}^{D}$ the sub-set
of such dominant weights. To make the notation homogenous, we write
$\omega_{p}^{B}=\Lambda_{p}^{B}$ for $p=1,...,n-1$ and $\omega_{p}^{D}%
=\Lambda_{p}^{D}$ for $p=1,...,n-2$.

\noindent Every dominant weight $\lambda=\underset{p=1}{\overset{n}{\sum}%
}\lambda_{p}\Lambda_{p}\in P_{+}$ has a unique decomposition
\begin{equation}
\lambda=\Lambda+\lambda^{\prime} \label{dec_lambda1}%
\end{equation}
such that $\lambda^{\prime}\in\Omega_{+}$ and
\begin{gather}
\left\{
\begin{tabular}
[c]{l}%
$\Lambda=\Lambda_{n}^{B}$ if $\lambda\in P_{+}^{B}$ and $\lambda_{n}$ is odd\\
$\Lambda=0$ if $\lambda\in P_{+}^{B}$ and $\lambda_{n}$ is even
\end{tabular}
\right.  \text{,}\label{decom_lambda}\\
\left\{
\begin{tabular}
[c]{l}%
$\Lambda=\Lambda_{n}^{D}\ $if $\lambda\in P_{+}^{D},$ $\lambda_{n}%
-\lambda_{n-1}\neq0$ $\operatorname{mod}2$ and $\lambda_{n}-\lambda_{n-1}>0$\\
$\Lambda=0$ if $\lambda\in P_{+}^{D}$ and $\lambda_{n}-\lambda_{n-1}=0$
$\operatorname{mod}2$\\
$\Lambda=\Lambda_{n-1}^{D}$ if $\lambda\in P_{+}^{D},$ $\lambda_{n}%
-\lambda_{n-1}\neq0$ $\operatorname{mod}2$ and $\lambda_{n}-\lambda_{n-1}<0$%
\end{tabular}
\right.  .\nonumber
\end{gather}

Now suppose $\lambda\in\Omega_{+}$.\ We write

\begin{description}
\item $Y_{\lambda}^{B}$ for the Young diagram having $\lambda_{i}$ columns of
height $i$ for $i=1,...,n-1$ and $\lambda_{n}/2$ columns of height $n$ if
$\lambda\in\Omega_{+}^{B}$

\item $Y_{\lambda}^{D}$ for the Young diagram having $\lambda_{i}$ columns of
height $i$ for $i=1,...,n-2$, $\min(\lambda_{n},\lambda_{n-1})$ columns of
height $n-1$ and $\left|  \lambda_{n}-\lambda_{n-1}\right|  /2$ columns of
height $n$ if $\lambda\in\Omega_{+}^{D}.$
\end{description}

\noindent If $\lambda\in\Omega_{+}^{D},$ $Y_{\lambda}^{D}$ may not suffice to
characterize the weight $\lambda$ because a column diagram of length $n$ may
be associated to $\omega_{n}^{D}$ or to $\overline{\omega}_{n}^{D}$. In the
sequel, we need to attach to $\lambda\in\Omega^{+}$ a combinatorial object
$Y(\lambda)$ analogous to the Young tableau associated to a dominant weight
for type $A$. This leads us to set:%
\begin{align}
\mathrm{(i)}  &  :Y^{B}(\lambda)=Y_{\lambda}^{B}\text{ if }\lambda\in
\Omega_{+}^{B},\nonumber\\
\mathrm{(ii)}  &  :Y^{D}(\lambda)=(Y_{\lambda}^{D},+)\text{ if }\lambda
\in\Omega_{+}^{D}\text{ and }\lambda_{n}-\lambda_{n-1}>0,\nonumber\\
\mathrm{(iii)}  &  :Y^{D}(\lambda)=(Y_{\lambda}^{D},0)\text{ if }\lambda
\in\Omega_{+}^{D}\text{ and }\lambda_{n}-\lambda_{n-1}=0,
\label{Def_Y(lambda)}\\
\mathrm{(iv)}  &  :Y^{D}(\lambda)=(Y_{\lambda}^{D},-)\text{ if }\lambda
\in\Omega_{+}^{D}\text{ and }\lambda_{n}-\lambda_{n-1}<0.\nonumber
\end{align}
Given $Y^{D}(\lambda)=(Y_{\lambda}^{D},\varepsilon)$ with $\varepsilon
\in\{+,0,-\}$ it is easy to recover $\lambda$.\ When $\varepsilon\in\{+,0\}$
(resp. $\varepsilon=-$), we have $\lambda=\underset{i=1}{\overset{n}{\sum}}%
\mu_{i}\omega_{i}$ (resp.\ $\lambda=\underset{i=1}{\overset{n-1}{\sum}}\mu
_{i}\omega_{i}+\mu_{n}\overline{\omega}_{n}$) where $\mu_{i}$ is the number of
columns of height $i$ in $Y_{\lambda}^{D}.$

\bigskip

Let us consider the ordered alphabets related to the crystal graphs
(\ref{vect_B}) and (\ref{vect_D})%
\begin{align*}
\mathcal{B}_{n}  &  =\{1\prec\cdot\cdot\cdot\prec n\prec0\prec\overline
{n}\prec\cdot\cdot\cdot\prec\overline{1}\}\text{ and }\\
\mathcal{D}_{n}  &  =\{1\prec\cdot\cdot\cdot\prec n-1\prec%
\begin{array}
[c]{l}%
n\\
\overline{n}%
\end{array}
\prec\overline{n-1}\prec\cdot\cdot\cdot\prec\overline{1}\}.
\end{align*}
Note that $\mathcal{D}_{n}$ is only partially ordered: $n$ and $\overline{n}$
are not comparable. We say that the letters of $\{1,...,n\}$ and
$\{\overline{1},...,\overline{n}\}$ are respectively unbarred and barred. We
set $\overline{0}=0$ and for any letter $x,$ $\overline{\overline{x}}=x.$
Write $\mathcal{B}_{n}^{\ast}$ and $\mathcal{D}_{n}^{\ast}$ for the free
monoids on $\mathcal{B}_{n}$ and $\mathcal{D}_{n}.$ Then we identify each
vertex $x_{1}\otimes\cdot\cdot\cdot\otimes x_{l}$ of $G_{n}$ with the
corresponding word $x_{1}\cdot\cdot\cdot x_{l}$ of $\mathcal{B}_{n}^{\ast}$ or
$\mathcal{D}_{n}^{\ast}.$

A column of type $B$ is a Young diagram
\[
C=%
\begin{tabular}
[c]{|l|}\hline
$x_{1}$\\\hline
$\cdot$\\\hline
$\cdot$\\\hline
$x_{l}$\\\hline
\end{tabular}
\]
of column shape filled by letters of $\mathcal{B}_{n}$ such that $C$ increases
from top to bottom and $0$ is the unique letter of $\mathcal{B}_{n}$ that may
appear more than once.

\noindent A column of type $D$ is a Young diagram $C$ of column shape filled
by letters of $\mathcal{D}_{n}$ such that $x_{i+1}\nleqslant x_{i}$ for
$i=1,...,l-1$. Note that the letters $n$ and $\overline{n}$ are the unique
letters that may appear more than once in $C$ and if they do, these letters
are different in two adjacent boxes. The height $h(C)$ of the column $C$ is
the number of its letters. The word obtained by reading the letters of $C$
from top to bottom is called the reading of $C$ and denoted by \textrm{w}%
$(C)$. Write $\mathbf{C}^{B}(n,p)$ (resp. $\mathbf{C}^{D}(n,p)$) for the set
of columns of height $p$ on $\mathcal{B}_{n}$ (resp.\ on $\mathcal{D}_{n})$.

Consider $\lambda^{\prime}\in\Omega^{B}.$ A tabloid of type $B$ and shape
$Y(\lambda^{\prime})$ is as a filling $\tau^{\prime}$ of the Young diagram
$Y_{\lambda^{\prime}}$ by letters of $\mathcal{B}_{n}$.\ Now consider
$\lambda^{\prime}\in\Omega^{D}.$ A tabloid of type $D$ and shape
$Y(\lambda^{\prime})$ is as a filling $\tau^{\prime}$ of the Young diagram
$Y_{\lambda^{\prime}}$ by letters of $\mathcal{D}_{n}$.\ A tabloid
$\tau^{\prime}$ of shape $\lambda^{\prime}$ can be regarded as the
juxtaposition $\tau^{\prime}=C_{1}\cdot\cdot\cdot C_{r}$ of its columns.\ The
reading of $\tau^{\prime}$ is $\mathrm{w(}\tau)=\mathrm{w(}C_{r})\cdot
\cdot\cdot\mathrm{w(}C_{1})\in G_{n}.$

\noindent Let $\lambda=\lambda^{\prime}+\Lambda\in P_{+}-\Omega_{+}$ be as in
(\ref{dec_lambda1})$.\;$A tabloid of shape $Y(\lambda)$ is a diagram obtained
by adding a spin column $\frak{C}$ in front of a tabloid of shape
$Y(\lambda^{\prime}).\;$The reading of the tabloid $\tau=\frak{C}C_{1}%
\cdot\cdot\cdot C_{r}$ of shape $Y(\lambda)$ is $\mathrm{w(}\tau
)=\frak{C}\mathrm{w(}C_{r})\cdot\cdot\cdot\mathrm{w(}C_{1})\in\frak{G}_{n}.$
For any dominant weight $\lambda\in P_{+}^{B}$ (resp. $P_{+}^{D}),$ we denote
by $\mathbf{T}^{B}(n,\lambda)$ (resp.\ $\mathbf{T}^{D}(n,\lambda))$ the set of
tabloids of shape $Y(\lambda).$

\begin{definition}
(Kashiwara-Nakashima)\label{defKN}

\begin{itemize}
\item  Let $\lambda\in P_{+}^{B}$.\ An orthogonal tableau of type $B$ and
shape $Y(\lambda)$ is a tabloid $T$ such that $\mathrm{w}(T)\in B(\lambda)$.

\item  Let $\lambda\in P_{+}^{D}$.\ An orthogonal tableau of type $D$ and
shape $Y(\lambda)$ is a tabloid $T$ such that $\mathrm{w}(T)\in B(\lambda)$.
\end{itemize}
\end{definition}

\noindent We will denote by $\mathbf{OT}^{B}(n,\lambda)$ and $\mathbf{OT}%
^{D}(n,\lambda)$ the sets of orthogonal tableaux respectively of type $B$ and
$D$ and shape $Y(\lambda).\;$An orthogonal tableau containing only one column
is called an admissible column. Set $\mathbf{Ca}^{B}(n,p)=\mathbf{OT}%
^{B}(n,\omega_{p}^{B})$ for $p\in\{1,...,n\},$ $\mathbf{Ca}^{D}%
(n,p)=\mathbf{OT}^{D}(n,\omega_{p}^{D})$ for $p\in\{1,...,n-1\},$
$\mathbf{Ca}_{+}^{D}(n,n)=\mathbf{OT}^{D}(n,\omega_{n}^{D})$ and
$\mathbf{Ca}_{-}^{D}(n,n)=\mathbf{OT}^{D}(n,\overline{\omega}_{n}^{D}).$ A
tabloid $T=C_{1}\cdot\cdot\cdot C_{r}$ of shape $Y(\lambda)$ with $\lambda
\in\Omega^{+}$ is an orthogonal tableau if and only if the tabloids of two
columns $C_{i}C_{i+1},$ $i=1,...,r-1$ are orthogonal tableaux \cite{KN}.\ This
is equivalent to say that its columns are admissible, its rows increase from
left to right for $\preceq$ and $T$ does not contain certain configurations
which depend on its type.

\bigskip

\noindent\textbf{Remark}: \label{rem_dupli}A column $C$ is admissible if and
only if it can be duplicated in a pair $(lC,rC)$ of columns which do not
contain any pair of letters $(z,\overline{z})$ following a certain process
described in \cite{L2(RS)}. This give an alternative combinatoric
caracterization of the orthogonal tableaux.\ The tabloid $T=C_{1}\cdot
\cdot\cdot C_{r}$ is an orthogonal tableau if and only if its columns can be
duplicated and the tabloid $\mathrm{spl}(T)=lC_{1}rC_{1}\cdot\cdot\cdot
lC_{r}rC_{r}$ is an orthogonal tableau.\ The advantage of this presentation is
that most of the forbidden configurations listed in \cite{KN} can not appear
in the duplicated form $\mathrm{spl}(T).$ In particulary, $T$ is an orthogonal
tableau of type $B$ if and only if the rows of $\mathrm{spl}(T)$ increase from
left to right for $\preceq.$

\noindent Consider $\lambda\in P_{+}-\Omega_{+}$ and $\frak{C}T$ a tabloid of
shape $Y(\lambda)$.\ Write $C_{\frak{C}}$ for the admissible column of height
$n$ containing the letters of $\frak{C}$.\ Then $\frak{C}T$ is an orthogonal
tableau if and only if $C_{\frak{C}}T$ is an orthogonal tableau.

\subsection{Canonical basis of a $U_{q}(\frak{g})$-module}

In the sequel we identify $B(\lambda)$ to $\{\mathrm{w}(T);$ $T\in
\mathbf{OT}(n,\lambda)\}$. Denote by $F\mapsto\overline{F}$ the involution of
$U_{q}(\frak{g})$ defined as the ring automorphism satisfying%
\[
\overline{q}=q^{-1},\text{ \ }t_{i}=t_{i}^{-1},\text{ \ \ }\overline{e_{i}%
}=e_{i},\text{ \ \ }\overline{f_{i}}=f_{i}\text{ \ \ \ for }i=1,...,n.
\]
By writing each vector $v$ of $V(\lambda)$ in the form $v=Fv_{\lambda}$ where
$F\in U_{q}(\frak{g})$, we obtain an involution of $V(\lambda)$ defined by%
\[
\overline{v}=\overline{F}v_{\lambda}.
\]
Let $U_{\mathbb{Q}}^{-}$ be the subalgebra of $U_{q}(\frak{g})$ generated over
$\mathbb{Q}[q,q^{-1}]$ by the $f_{i}^{(k)}$ and set $V_{\mathbb{Q}}%
(\lambda)=U_{\mathbb{Q}}^{-}v_{\lambda}$. We can now state:

\begin{theorem}
\label{TH_K_2}(Kashiwara) There exists a unique $\mathbb{Q[}q,q^{-1}]$-basis
$\{G(T);$ $T\in\mathbf{OT}(n,\lambda)\}$ of $V_{\mathbb{Q}}(\lambda)$ such
that:%
\begin{gather}
G(T)\equiv\mathrm{w(}T)\text{ }\mathrm{mod}\text{ }qL(\lambda
),\label{cond_cong}\\
\overline{G(T)}=G(T). \label{cond_invo}%
\end{gather}
\end{theorem}

\noindent Note that $G(T)\in V_{\mathbb{Q}}(\lambda)\cap L(\lambda)$.\ The
basis $\{G(T);T\in\mathbf{OT}(n,\lambda)\}$ is called the lower global (or
canonical) basis of $V(\lambda)$, and our aim is to calculate it.

\subsection{Marsh's algorithm\label{sub-sec_M_algo}}

Now we review Marsh's algorithm \cite{Ma} for computing the global basis of
$V(\omega_{p})$. Let $\mathrm{w(}C)\in B(\omega_{p}).$ The letter
$x\in\mathrm{w(}C)$ is movable in $\mathrm{w(}C)$ for $i$ if $\widetilde
{e}_{i}(x)\neq0$ and is not a letter of $\mathrm{w(}C)$. We define a path in
$B(\omega_{p})$ joining $\mathrm{w(}C)$ to $b_{\omega_{p}}$.

\noindent If $\mathrm{w}(C)\neq b_{\omega_{p}},$ let $z$ be the leftmost
movable letter of $\mathrm{w(}C).\;$When $C$ is of type $B,$ let $i_{1}$ be
the unique $i\in\{1,...,n\}$ such that $\widetilde{e}_{i}(z)\neq0.$ When $C$
is of type $D$ we define $i_{1}$ as follows:%
\[%
\begin{tabular}
[c]{l}%
when $z=\overline{n-1}$, $i_{1}=n-1,$\\
when $z=\overline{n},$ $i_{1}=\left\{
\begin{array}
[c]{l}%
n-1\text{ if }\mathrm{w}_{n-1}(C)=(\overline{n}n)^{r}(\overline{n-1})\\
n\text{ otherwise}%
\end{array}
\right.  ,$\\
when $z=n,$ $i_{1}=\left\{
\begin{array}
[c]{l}%
n\text{ if }\mathrm{w}_{n}(C)=(n\overline{n})^{r}(\overline{n-1})\\
n-1\text{ otherwise}%
\end{array}
\right.  ,$\\
$i_{1}$ the unique $i$ such that $\widetilde{e}_{i}(z)\neq0$ in the other
cases.
\end{tabular}
\]
Note that the admissibility of $C$ and the definition of $z$ imply that the
letters $n$ and $\overline{n}$ do not belong to $C$ when $z=\overline{n-1}.$
Marsh has proved that $i_{1}$ is always such that $\widetilde{e}_{i_{1}%
}(\mathrm{w}(C))\neq0$. Set $p_{1}=\varepsilon_{i_{1}}(\mathrm{w}%
(C))\in\{1,2\}$. Let $C_{1}$ be the admissible column of reading
$\widetilde{e}_{i_{1}}^{p_{1}}(\mathrm{w}(C)).$ Then we compute similarly
$C_{2}$ from $C_{1}$ and write $\mathrm{w(}C_{2})=\widetilde{e}_{i_{2}}%
^{p_{2}}\mathrm{w(}C_{1})$. Finally, after a finite number of steps, we will
reach $b_{\omega_{p}}$ and we will get $\mathrm{w(}b_{\omega_{p}}%
)=\widetilde{e}_{i_{r}}^{p_{r}}\cdot\cdot\cdot\widetilde{e}_{i_{1}}^{p_{1}%
}\mathrm{w(}C)$, hence $\mathrm{w(}C)=\widetilde{f}_{i_{1}}^{p_{1}}\cdot
\cdot\cdot\widetilde{f}_{i_{r}}^{p_{r}}b_{\omega_{p}}$.

In fact $\widetilde{f}_{i_{1}}^{p_{1}},...,\widetilde{f}_{i_{r}}^{p_{r}}$ are
chosen to verify:%
\begin{equation}
f_{i_{1}}^{(p_{1})}\cdot\cdot\cdot f_{i_{r}}^{(p_{r})}v_{\Lambda_{p}%
}=\widetilde{f}_{i_{1}}^{p_{1}}\cdot\cdot\cdot\widetilde{f}_{i_{r}}^{p_{r}%
}v_{\Lambda_{p}}. \label{eq_G_C_for_fund}%
\end{equation}
where for any integer $m$, $e_{i}^{(m)}=e_{i}^{m}/[m]_{i}!$ and $f_{i}%
^{(m)}=f_{i}^{m}/[m]_{i}!$. This implies

\begin{theorem}
\label{Th_marsh}(Marsh) For any admissible column $C$%
\[
G(C)=f_{i_{1}}^{(p_{1})}\cdot\cdot\cdot f_{i_{r}}^{(p_{r})}v_{\omega_{p}%
}\text{,}%
\]
where the integers $i_{1},...,i_{r}$ and $p_{1},...,p_{r}$ are determined by
the algorithm above.
\end{theorem}

\section{Wedge products of the vector representation}

In this section we introduce a $q$-analogue $W(\omega_{p})$ of the $p$-th
wedge product of the vector representation of $U_{q}(\frak{g})$.\ It has been
inspired by a similar construction for quantized affine Lie algebras obtained
by Kashiwara, Miwa, Petersen and Yung in \cite{KMPY}. When $q\rightarrow1$ in
$W(\omega_{p}),$ we do not recover the classical notion of the $p$-th wedge
product. In particular $W(\omega_{p})$ is not irreducible.\ We are going to
see that it contains an irreducible component isomorphic to $V(\omega_{p}).$

\noindent The vector representation $V(\Lambda_{1}^{B})$ of $U_{q}%
(\frak{so}_{2n+1})$ is the vector space of basis $\{v_{x},$ $x\in
\mathcal{B}_{n}\}$ where
\[
t_{i}(v_{x})=q_{i}^{<\mathrm{wt}(x),\alpha_{i}>}\text{ for }i=1,...,n,
\]%
\[
\left\{
\begin{tabular}
[c]{l}%
$e_{n}(v_{\overline{n}})=v_{0},$ $e_{n}(v_{0})=(q+q^{-1})v_{n}$ and
$e_{n}(v_{x})=0$ if $x\notin\{0,n\}$\\
$f_{n}(v_{n})=v_{0},\text{ }f_{n}(v_{0})=(q+q^{-1})v_{\overline{n}}\text{ and
}e_{n}(v_{x})=0\text{ if }x\notin\{0,\overline{n}\}$%
\end{tabular}
\right.
\]
and the action of $e_{i}$, $f_{i}$, $i=1,...,n-1$ is determined by that of
$\widetilde{e}_{i}$, $\widetilde{f}_{i}$ on the crystal $B(\Lambda_{1}^{B})$,
that is%
\begin{equation}
\left\{
\begin{tabular}
[c]{l}%
$f_{i}(v_{x})=v_{y}$ if $\widetilde{f}_{i}(x)=y$, $f_{i}(v_{x})=0$ otherwise\\
$e_{i}(v_{x})=v_{y}$ if $\widetilde{e}_{i}(x)=y$, $e_{i}(v_{x})=0$ otherwise
\end{tabular}
\right.  . \label{action fi_vect}%
\end{equation}
Note that, with our definition of the action of $f_{n}$ on $v_{0},$ we have
$f_{n}^{(2)}(v_{n})=v_{\overline{n}}$.

\noindent The vector representation $V(\Lambda_{1}^{D}$) of $U_{q}%
(\frak{so}_{2n+1})$ is the space vector of basis $\{v_{x},$ $x\in
\mathcal{D}_{n}\}$ where the action of $e_{i},$ $f_{i}$ and $t_{i}$ ,
$i=1,...,n$ is given by $t_{i}(v_{x})=q_{i}^{<\mathrm{wt}(x),\alpha_{i}>}$ and
(\ref{action fi_vect}).

\subsection{The representations $W(\omega_{p}^{B})$ and $W(\omega_{p}^{D})$}

We start with the decomposition of the square tensor product of $V(\Lambda
_{1})$%
\begin{equation}
V(\Lambda_{1})^{\otimes2}\cong M%
%TCIMACRO{\tbigoplus }%
%BeginExpansion
{\textstyle\bigoplus}
%EndExpansion
N%
%TCIMACRO{\tbigoplus }%
%BeginExpansion
{\textstyle\bigoplus}
%EndExpansion
V(0) \label{dec_square_of B(1)}%
\end{equation}
where $N\cong V(2\Lambda_{1})$ the submodule generated by $v_{1}\otimes
v_{1},$ $V(0)$ is the trivial representation and
\[
M=\left\{
\begin{tabular}
[c]{l}%
$V(\omega_{2}^{B})$ for type $B$ and $n=2$\\
$V(\omega_{2}^{D})%
%TCIMACRO{\tbigoplus }%
%BeginExpansion
{\textstyle\bigoplus}
%EndExpansion
V(\overline{\omega}_{2}^{D})$ for type $D$ and $n=2$\\
$V(\omega_{2}^{D})$ for type $D$ and $n=3$\\
$V(\Lambda_{2})$ otherwise
\end{tabular}
\right.  .
\]
We set%
\begin{equation}
W(\omega_{p}^{B})=V(\Lambda_{1}^{B})^{\otimes p}/N_{p}^{B}\text{ and }%
W(\omega_{p}^{D})=V(\Lambda_{1}^{D})^{\otimes p}/N_{p}^{D} \label{def_q_wedge}%
\end{equation}
where%

\[
N_{p}=\underset{i=0}{\overset{p-2}{\sum}}V(\Lambda_{1})^{\otimes i}\otimes
N\otimes V(\Lambda_{1})^{\otimes(p-2-i)}.
\]
The representations $W(\omega_{p}^{B})$ and $W(\omega_{p}^{D})$ are
respectively called the $p$-th $q$-wedge power of $V(\Lambda_{1}^{B})$ and
$V(\Lambda_{1}^{D}).$ Let $\Psi_{p}$ be the canonical projection
$V(\Lambda_{1})^{\otimes p}\rightarrow W(\omega_{p})$.\ Denote by $v_{x_{1}%
}\wedge\cdot\cdot\cdot\wedge v_{x_{p}}$ the image of the vector $v_{x_{1}%
}\otimes\cdot\cdot\cdot\otimes v_{x_{p}}$ by $\Psi_{p}$.

\begin{proposition}
\ \ \ \ \ \ \ \ \ \ \label{prop_relation}

\noindent In $W(\omega_{2}^{B})$ we have the relations

\begin{enumerate}
\item $v_{x}\wedge v_{x}=0$ for $x\neq0$,

\item $v_{y}\wedge v_{x}=-q^{2}v_{x}\wedge v_{y}$ for $x\neq\overline{y}$ and
$x\prec y,$

\item $v_{\overline{i}}\wedge v_{i}=-q^{4}v_{i}\wedge v_{\overline{i}%
}+(1-q^{4})\underset{k=1}{\overset{n-i}{\sum}}(-1)^{k}q^{2k}v_{i+k}\wedge
v_{\overline{i+k}}+(-1)^{n-i+1}q^{2(n-i)+1}v_{0}\wedge v_{0}$ for $i=1,...,n.$
\end{enumerate}

\noindent In $W(\omega_{2}^{D})$ we have the relations

\begin{enumerate}
\item $v_{x}\wedge v_{x}=0$,

\item $v_{y}\wedge v_{x}=-qv_{x}\wedge v_{y}$ for $x\neq\overline{y}$ and
$x\prec y,$

\item $v_{\overline{i}}\wedge v_{i}=-q^{2}v_{i}\wedge v_{\overline{i}%
}+(1-q^{2})\underset{k=1}{\overset{n-i-1}{\sum}}(-q)^{k}v_{i+k}\wedge
v_{\overline{i+k}}+(-q)^{n-i}(v_{n}\wedge v_{\overline{n}}+v_{\overline{n}%
}\wedge v_{n})$ for $i=1,...,n-1.$
\end{enumerate}
\end{proposition}

\begin{proof}
Using (\ref{tensor1}), (\ref{tensor2}) and (\ref{tensor3}) we obtain that the
following vectors belong to $N^{B}:$

\noindent$v_{x}\otimes v_{x}$ with $x\in\mathcal{B}_{n}-\{0\},$

\noindent$v_{\overline{n}}\otimes v_{n}+q^{4}v_{n}\otimes v_{\overline{n}%
}+qv_{0}\otimes v_{0},$

\noindent$v_{y}\otimes v_{x}+q^{2}v_{x}\otimes v_{y}$ with $x,y\in
\mathcal{B}_{n}$ such that $x\prec y$ and $x\neq\overline{y},$

\noindent$v_{\overline{k}}\otimes v_{k}+q^{2}v_{\overline{k+1}}\otimes
v_{k+1}+q^{2}v_{k+1}\otimes v_{\overline{k+1}}+q^{4}v_{k}\otimes
v_{\overline{k}}$ with $k\in\{1,...,n-1\}.$

\noindent For example $v_{i+1}\otimes v_{i+1}=f_{i}^{(2)}\left(  v_{i}\otimes
v_{i}\right)  $ for $i=1,...,n-1,$ $f_{n}^{(2)}\left(  v_{n}\otimes
v_{n}\right)  =v_{\overline{n}}\otimes v_{n}+q^{4}v_{n}\otimes v_{\overline
{n}}+qv_{0}\otimes v_{0},$ $f_{n}^{(4)}\left(  v_{n}\otimes v_{n}\right)
=v_{\overline{n}}\otimes v_{\overline{n}}$ and $v_{\overline{i}}\otimes
v_{\overline{i}}=f_{i}^{(2)}\left(  v_{\overline{i+1}}\otimes v_{\overline
{i+1}}\right)  $ for $i=1,...,n-1$.

\noindent Hence the relations

\noindent$\mathrm{(i):}$ $v_{x}\wedge v_{x}=0$ with $x\in\mathcal{B}_{n}-\{0\},$

\noindent$\mathrm{(ii):}$ $v_{\overline{n}}\wedge v_{n}+q^{4}v_{n}\wedge
v_{\overline{n}}+qv_{0}\wedge v_{0}=0,$

\noindent$\mathrm{(iii):}$ $v_{y}\wedge v_{x}+q^{2}v_{x}\wedge v_{y}=0$ with
$x\in\mathcal{B}_{n}$ such that $x\neq\overline{y}$ and $x\prec y,$

\noindent$\mathrm{(iv):}$ $v_{\overline{k}}\wedge v_{k}+q^{2}v_{\overline
{k+1}}\wedge v_{k+1}+q^{2}v_{k+1}\wedge v_{\overline{k+1}}+q^{4}v_{k}\wedge
v_{\overline{k}}=0$ with $k\in\{1,...,n-1\},$

\noindent are true in $W(\omega_{2}^{B}).$ Then Relation $3$ in $W(\omega
_{2}^{B})$ may be obtained by a straightforward induction from relations
$\mathrm{(iii)}$ and $\mathrm{(iv)}$ below.

To obtain the relations $1,2$ and $3$ in $W(\omega_{2}^{D}),$ we consider
similarly the vectors

\noindent$v_{x}\otimes v_{x}$ with $x\in\mathcal{D}_{n},$

\noindent$v_{y}\otimes v_{x}+qv_{x}\otimes v_{y}$ with $x,y\in\mathcal{D}_{n}$
such that $x\prec y$ and $x\neq\overline{y},$

\noindent$v_{\overline{k}}\otimes v_{k}+qv_{\overline{k+1}}\otimes
v_{k+1}+qv_{k+1}\otimes v_{\overline{k+1}}+q^{2}v_{k}\otimes v_{\overline{k}}$
with $k\in\{1,...,n-1\},$

\noindent that belong to $N^{D}.$
\end{proof}

For any column $C$ of reading $w=c_{1}\cdot\cdot\cdot c_{p}$ where the $c_{i}%
$'s are letters, we set $v_{C}=v_{c_{1}}\wedge\cdot\cdot\cdot\wedge v_{c_{p}}%
$.\ Then each vector $\Psi_{p}(v_{x_{1}}\otimes\cdot\cdot\cdot\otimes
v_{x_{p}})=v_{x_{1}}\wedge\cdot\cdot\cdot\wedge v_{x_{p}}$ can be decomposed
into a linear combination of vectors $v_{C}$ by applying from left to right a
sequence of relations given in the above proposition.\ 

\begin{lemma}
The vectors of $\{v_{C},$ $C\in C(n,p)\}$ form a basis of $W(\omega_{p}).$
\end{lemma}

\begin{proof}
We only sketch the proof for type $B.\;$The arguments are essentially the same
for type $D$.\ Each vector of $W(\omega_{p}^{B})$ can be decomposed into a
linear combination of vectors $v_{C}.\;$Hence $\dim(W(\omega_{p}^{B}%
))\leq\mathrm{card}(C(n,p)).$ For any integer $k$ such that $0\leq k\leq p/2,$
there is $\binom{2n+1}{p-2k}$ columns of type $B$ containing $2k$ or $2k+1$
letters $0.$ So $\mathrm{card}(C(n,p))=\underset{0\leq k\leq p/2}{\sum}%
\binom{2n+1}{p-2k}$ and $\dim(W(\omega_{p}^{B}))\leq\underset{0\leq k\leq
p/2}{\sum}\binom{2n+1}{p-2k}.$

Accordingly to (\ref{dec_square_of B(1)}) there exists in $V(\Lambda
_{1})^{\otimes2}$ a highest weight vector $v_{\emptyset}$ of weight $0$.\ For
any integer $k$ such that $0\leq k\leq p/2,$ set%
\[
u_{k}=v_{1}\otimes\cdot\cdot\cdot\otimes v_{p-2k}\otimes v_{\emptyset
}^{\otimes k}\subset V(\Lambda_{1})^{\otimes p}.
\]
Then $\Psi_{p}(u_{k})$ is of highest weight $\omega_{p-2k}^{B}$ (with
$\omega_{0}^{B}=0$ if $p=2k$) in $W(\omega_{p}^{D}).\;$So for any integer $k$
such that $0\leq k\leq p/2,$ $W(\omega_{p}^{B})$ contains at least an
irreducible component isomorphic to $V(\omega_{p-2k}^{B})$.\ Hence
$\dim(W(\omega_{p}^{B}))\geq\underset{0\leq k\leq p/2}{\sum}\binom{2n+1}%
{p-2k}$ since $\dim(V(\omega_{p-2k}^{B}))=\binom{2n+1}{p-2k}.$ Finally
$\dim(W(\omega_{p}^{B}))=\mathrm{card}(C(n,p))$ and $\{v_{C},$ $C\in C(n,p)\}$
is a basis of $W(\omega_{p}^{B}).$ Note that we have the decomposition%
\[
W(\omega_{p})\cong%
%TCIMACRO{\dbigoplus _{0\leq k\leq p/2}}%
%BeginExpansion
{\displaystyle\bigoplus_{0\leq k\leq p/2}}
%EndExpansion
V(\omega_{p-2k}).
\]
\end{proof}

The coordinates of a vector $v_{x_{1}}\wedge\cdot\cdot\cdot\wedge v_{x_{p}}$
on the basis $\{v_{C},$ $C\in\mathbf{C}(n,p)\}$ are all in $\mathbb{Z[}q]$
since it is true for the coefficients appearing in the relations of
Proposition \ref{prop_relation}.

\noindent Consider the $A$-lattice $L_{p}=\underset{C\in\mathbf{C}(n,p)}{%
%TCIMACRO{\tbigoplus }%
%BeginExpansion
{\textstyle\bigoplus}
%EndExpansion
}Av_{C}$ of $W(\omega_{p})$ and denote by $\pi_{p}$ the projection
$L_{p}\rightarrow L_{p}/qL_{p}$.\ We identify $\pi_{p}(v_{C})$ with the word
$\mathrm{w}(C)$.

\begin{lemma}
\label{lem_base_cryst_Wp}$(L_{p},B_{p}=\{\mathrm{w}(C),$ $C\in\mathbf{C}%
(n,p)\})$ is a crystal basis of $W(\omega_{p})$.
\end{lemma}

\begin{proof}
Note first that $L_{p}=\underset{\mu\in P}{%
%TCIMACRO{\tbigoplus }%
%BeginExpansion
{\textstyle\bigoplus}
%EndExpansion
}L_{\mu}$ and $B_{p}=\underset{\mu\in P}{\cup}B_{\mu}$.\ For any word
$w=x_{1}\cdot\cdot\cdot x_{p}$ of length $p$ set $v_{w}=v_{x_{1}}\otimes
\cdot\cdot\cdot\otimes v_{x_{p}}.\;$Consider $\mathcal{L}_{p}=\underset{w}{%
%TCIMACRO{\tbigoplus }%
%BeginExpansion
{\textstyle\bigoplus}
%EndExpansion
}Av_{w}$ where the sum runs over the set of words of length $p.$ Then
$(\mathcal{L}_{p},\mathcal{B}_{p}=\{w\in G_{n},$ $l(w)=p\})$ is a crystal
basis of $V(\Lambda_{1})^{\otimes p}$ because it is the $p$-th tensor power of
the crystal basis of $V(\Lambda_{1})$.

One has $\Psi_{p}(\mathcal{L}_{p})=L_{p}.$ Indeed, in the decomposition of
$\Psi_{p}(v_{x_{1}}\otimes\cdot\cdot\cdot\otimes v_{x_{p}})=v_{x_{1}}%
\wedge\cdot\cdot\cdot\wedge v_{x_{p}}$ on the basis $\{v_{C},$ $C\in
\mathbf{C}(n,p)\}$ the coefficients are all in $\mathbb{Z}[q]\subset A$.\ For
any $v\in\mathcal{L}_{p}$%
\[
\widetilde{f}_{i}(\Psi_{p}(v))=\Psi_{p}(\widetilde{f}_{i}(v))\text{ and
}\widetilde{e}_{i}(\Psi_{p}(v))=\Psi_{p}(\widetilde{e}_{i}(v))
\]
because $N_{p}$ is a sub-module of $V(\Lambda_{1})^{\otimes p}$ so is stable
under the action of Kashiwara's operators. This implies that $\widetilde
{e}_{i}L_{p}\subset L_{p}$ and $\widetilde{f}_{i}L_{p}\subset L_{p}$ for
$i=1,...,n$. It is clear that $B_{p}$ is a $\mathbb{Q}$-basis of $L_{p}%
/qL_{p}$.\ The action of Kashiwara's operators on $B_{p}$ is obtained by
restriction of the action of these operators on $\mathcal{B}_{p}$.\ We deduce
from (\ref{TENS1}) and (\ref{TENS2}) that the image of a column word
containing $p$ letters by any $\widetilde{f}_{i}$ or $\widetilde{e}_{i},$
$i=1,...,n$ is zero or a column word containing $p$ letters.\ So
$\widetilde{e}_{i}B\subset B\cup\{0\}$, $\widetilde{f}_{i}B\subset B\cup
\{0\}$.\ Moreover for any $b_{1},b_{2}\in B_{p}$ and $i\in\{1,...,n\},$
$\widetilde{e}_{i}b_{1}=b_{2}\Longleftrightarrow b_{1}=\widetilde{f}_{i}%
b_{2}.\;$This proves that $(L_{p},B_{p})$ is a crystal basis of $W(\omega
_{p})$ (see Definition 4.1 of \cite{Ka2})$.$
\end{proof}

For any $p=1,...,n,$ the vector $v_{\omega_{p}^{B}}=v_{1}\wedge\cdot\cdot
\cdot\wedge v_{p}$ is of highest weight $\omega_{p}^{B}$ in $W(\omega_{p}%
^{B})$.\ Similarly the vector $v_{\omega_{p}^{D}}=v_{1}\wedge\cdot\cdot
\cdot\wedge v_{p}$ is of highest weight $\omega_{p}^{D}$ in $W(\omega_{p}%
^{D})$ and the vector $v_{\overline{\omega}_{n}^{D}}=v_{1}\wedge\cdot
\cdot\cdot\wedge v_{\overline{n}}$ is of highest weight $\overline{\omega}%
_{n}^{D}$ in $W(\omega_{n}^{D}).$ For $p=1,...,n,$ we identify $V(\omega_{p})$
with the sub-module of $W(\omega_{p})$ generated by $v_{\omega_{p}}$ and
$V(\overline{\omega}_{n}^{D})$ with the sub-module of $W(\omega_{n}^{D})$
generated by $v_{\overline{\omega}_{n}^{D}}.$

\bigskip

\noindent\textbf{Remark}:

\begin{enumerate}
\item  Set $\mathcal{N}^{B}=N^{B}\bigoplus V(0)$ and $\mathcal{N}^{D}%
=N^{D}\bigoplus V(0)$.\ If we replace $N^{B}$ and $N^{D}$ by $\mathcal{N}^{B}$
and $\mathcal{N}^{D}$ in the definition of $W(\omega_{p}^{B})$ and
$W(\omega_{p}^{D}),$ we obtain $\mathcal{W}(\omega_{p}^{B})$ and
$\mathcal{W}(\omega_{p}^{D})$ the $q$-analogue of the $p$-th wedge power
introduced in \cite{JMO}.\ The representations $\mathcal{W}(\omega_{p}^{B})$
and $\mathcal{W}(\omega_{p}^{D})$ are irreducible and when $q\rightarrow1,$ we
recover the classical notion of the $p$-th wedge product of the vector
representation.\ Moreover the relations of Proposition \ref{prop_relation}
together with the relation $v_{1}\wedge v_{\overline{1}}=-v_{\overline{1}%
}\wedge v_{1}$ are verified in $\mathcal{W}(\omega_{p}^{B})$ and
$\mathcal{W}(\omega_{p}^{D})$. It is possible to show that every vector
$v_{x_{1}}\wedge\cdot\cdot\cdot\wedge v_{x_{p}}$ decomposes on the basis
$\{v_{C},$ $C\in\mathbf{Ca}(n,p)\}.\;$But this time the coefficients appearing
during this decomposition may not be in $\mathbb{Z}[q].$ For example when
$n=2$ the column
\begin{tabular}
[c]{|l|}\hline
$\mathtt{1}$\\\hline
$\mathtt{\bar{1}}$\\\hline
\end{tabular}
is not admissible and one has in $\mathcal{W}(\omega_{p})$%
\begin{align}
v_{1}\wedge v_{\overline{1}}  &  =q^{2}v_{2}\wedge v_{\overline{2}}%
+\frac{q^{3}}{q^{4}-1}v_{0}\wedge v_{0}\text{ for type }B,\label{1_1bar_for_B}%
\\
v_{1}\wedge v_{\overline{1}}  &  =\frac{q}{1-q^{2}}(v_{2}\wedge v_{\overline
{2}}+v\overline{_{2}}\wedge v_{2})\text{ for type }D. \label{1_1bar_for_D}%
\end{align}

Thus the modules $\mathcal{W}(\omega_{p})$ are certainly less convenient than
the $W(\omega_{p})$ for calculating the canonical basis, because the analog of
Lemma \ref{Lem_F()} below does not hold for $\mathcal{W}(\omega_{p}).\;$See
remarks below.

\item  It is possible to define for any $p=1,...,n$ the $p$-th $q$-wedge
product $W(\Lambda_{p}^{C})$ of the vector representation of $U_{q}(sp_{2n})$
in a similar way.\ The representation $W(\Lambda_{p}^{C})$ contains an
irreducible component isomorphic to the fundamental representation
$V(\Lambda_{p}^{C})$ and has a natural basis indexed by the columns of type
$C$ (which are the columns of type $B$ filled only by letters of
$\mathcal{B}_{n}-\{0\}$) and height $p.$ This module $W(\Lambda_{p}^{D})$ is
in fact that used in \cite{L1}. There exists in $V(\Lambda_{p}^{C})$ relations
analogous to those given in Proposition \ref{prop_relation}.\ They are
obtained from the relations holding in $V(\omega_{p}^{B})$ by setting
$v_{0}\wedge v_{0}=0.$
\end{enumerate}

\subsection{Action of the Chevalley's generators $f_{i}$ and $t_{i}$ on
$W(\omega_{p}).\label{sub_sec_formulas}$}

In the sequel we will need of a combinatorial description of the action of the
Chevalley generators $f_{i}$ and $t_{i},$ $i=1,...,n$ on the basis $\{v_{C},$
$C\in C(n,p)\}.$ The action of $t_{i},$ $i=1,...,n$ is given by%

\begin{equation}
t_{i}v_{C}=q_{i}^{<h_{i},v_{C}>}\,v_{C}. \label{hi}%
\end{equation}
Given a column $C,$ we write $\mathrm{w}_{i}(C)$ for the column word obtained
by reading only the letters of $C$ that occur in

\begin{description}
\item $\{\overline{i},\overline{i+1},i,i+1\}$ if $i=1,...,n-1,$

\item $\{\overline{n},0,n\}$ if $i=n$ and $C$ is of type $B,$

\item $\{\overline{n-1},\overline{n},n-1,n\}$ if $i=n$ and $C\,$is of type $D.$
\end{description}

\noindent If the vector $v_{Y}$ occur in the decomposition of $f_{i}(v_{C})$
on the basis $\{v_{X},$ $X\in C(n,p)\}$, the words $\mathrm{w}(C)$ and
$\mathrm{w}(Y)$ can only differ by letters of $\mathrm{w}_{i}(C)$ and
$\mathrm{w}_{i}(Y)$.\ Hence $Y$ is characterized by $\mathrm{w}_{i}(Y)$ and
the letters of $C$ which are not in $\mathrm{w}_{i}(C).$

\noindent To obtain the action of the operators $f_{i},$ it suffices to use
(\ref{tensor3}) and the relations obtained in Proposition \ref{prop_relation}.

\subsubsection{Action of the $f_{i}$'s on $W(\omega_{p}^{B})$}

Consider a column $C$ of type $B$. Suppose first $i\in\{1,...,n-1\}$.\ One
has
\begin{equation}
f_{i}(v_{C})=\left\{
\begin{tabular}
[c]{l}%
$q_{i}^{-1}v_{X}$ with $\mathrm{w}_{i}(X)=(i+1)\overline{i}$ if $\mathrm{w}%
_{i}(C)=(i+1)(\overline{i+1})$\\
\\
$v_{X}$ with $\mathrm{w}_{i}(X)=(i+1)\overline{i}$ if $\mathrm{w}%
_{i}(C)=i\overline{i}$\\
\\
$v_{X_{1}}+q_{i}v_{X_{2}}$ with$\left\{
\begin{tabular}
[c]{l}%
$\mathrm{w}_{i}(X_{1})=(i+1)(\overline{i+1})$\\
$\mathrm{w}_{i}(X_{2})=i\overline{i}$%
\end{tabular}
\right.  $if $\mathrm{w}_{i}(X)=i(\overline{i+1})$\\
\\
$v_{X}$ with $\mathrm{w}(X)=\widetilde{f}_{i}(\mathrm{w}(C))$ if $\varphi
_{i}(\mathrm{w}(C))=1$ and $\mathrm{w}_{i}(C)\neq(i+1)(\overline{i+1})$\\
\\
$0$ otherwise.
\end{tabular}
\right.  . \label{fi_B}%
\end{equation}
If $i=n,$ set $E_{n}=\{n,0,\overline{n}\}$. One has
\begin{equation}
f_{n}(v_{C})=\left\{
\begin{tabular}
[c]{l}%
$\frac{(1-(-q^{2})^{r})}{q}v_{X}$ with $\mathrm{w}_{n}(X)=0^{r-1}\overline{n}$
if $\mathrm{w}_{n}(C)=0^{r}$ and $r\geq1$\\
\\
$v_{X_{1}}+q(1-(-q^{2})^{r})v_{X_{2}}$ with$\left\{
\begin{tabular}
[c]{l}%
$\mathrm{w}_{n}(X_{1})=0^{r+1}$\\
$\mathrm{w}_{n}(X_{2})=n0^{r-1}\overline{n}$%
\end{tabular}
\right.  $if $\mathrm{w}_{n}(C)=n0^{r}$ and $r\geq1$\\
\\
$v_{X}$ with $\mathrm{w}_{n}(X)=0$ if $\mathrm{w}_{n}(X)=n$\\
\\
$v_{X}$ with $\mathrm{w}_{n}(X)=0^{r+1}\overline{n}$ if $\mathrm{w}%
_{n}(C)=n0^{r}\overline{n}$ and $r\geq0$\\
\\
$0$ otherwise
\end{tabular}
\right.  . \label{fn_B}%
\end{equation}
where for any $p\in\mathbb{N}$, $0^{p}$ means $0$ repeated $p$ times.

\subsubsection{Action of the $f_{i}$'s on $W(\omega_{p}^{D})$}

When $i\neq n-1,n,$ the action of $f_{i}$ on $v_{C}$ is given by (\ref{fi_B})
for any column $C$ of type $D$. For $i=n-1$, we obtain
\begin{equation}
f_{n-1}(v_{C})=\left\{
\begin{tabular}
[c]{l}%
$-q^{2r-1}v_{X}$ with $\mathrm{w}_{n-1}(X)=n(\overline{n}n)^{r-1}%
(\overline{n-1})$ if $\mathrm{w}_{n-1}(C)=(\overline{n}n)^{r}$ and $r\geq1$\\
\\
$v_{X_{1}}-q^{2r}v_{X_{2}}$ with $\left\{
\begin{tabular}
[c]{l}%
$\mathrm{w}_{n-1}(X_{1})=n(\overline{n}n)^{r}$\\
$\mathrm{w}_{n-1}(X_{2})=\left(  n-1\right)  n(\overline{n}n)^{r-1}%
(\overline{n-1})$%
\end{tabular}
\right.  $\\
\ \ \ \ \ \ \ \ \ \ \ \ \ \ \ \ \ \ \ \ \ \ \ \ \ \ \ \ \ \ \ \ \ \ \ \ \ \ \ \ \ \ \ \ \ \ \ \ \ \ \ \ \ \ \ \ \ \ \ \ \ \ if
$\mathrm{w}_{n-1}(C)=(n-1)(\overline{n}n)^{r}$ and $r\geq1$\\
\\
$v_{X_{1}}+q^{2r}v_{X_{2}}$ with $\left\{
\begin{tabular}
[c]{l}%
$\mathrm{w}_{n-1}(X_{1})=(\overline{n}n)^{r}(\overline{n-1})$\\
$\mathrm{w}_{n-1}(X_{2})=(n\overline{n})^{r}(\overline{n-1})$%
\end{tabular}
\right.  $if $\mathrm{w}_{n-1}(C)=(\overline{n}n)^{r}\overline{n}$ and
$r\geq1$\\
\\
$v_{X_{1}}+qv_{X_{2}}+q^{2r+1}v_{X_{3}}$ with $\left\{
\begin{tabular}
[c]{l}%
$\mathrm{w}_{n-1}(X_{1})=(n\overline{n})^{r+1}$\\
$\mathrm{w}_{n-1}(X_{2})=\left(  n-1\right)  (\overline{n}n)^{r}%
(\overline{n-1})$\\
$\mathrm{w}_{n-1}(X_{3})=\left(  n-1\right)  (n\overline{n})^{r}%
(\overline{n-1})$%
\end{tabular}
\right.  $\\
\ \ \ \ \ \ \ \ \ \ \ \ \ \ \ \ \ \ \ \ \ \ \ \ \ \ \ \ \ \ \ \ \ \ \ \ \ \ \ \ \ \ \ \ \ \ \ \ \ \ \ \ \ \ \ \ \ \ \ \ \ \ if
$\mathrm{w}_{n-1}(C)=(n-1)(\overline{n}n)^{r}\overline{n}$ and $r\geq1$\\
\\
$v_{X_{1}}+qv_{X_{2}}$ with $\left\{
\begin{tabular}
[c]{l}%
$\mathrm{w}_{n-1}(X_{1})=n\overline{n}$\\
$\mathrm{w}_{n-1}(X_{2})=\left(  n-1\right)  (\overline{n-1})$%
\end{tabular}
\right.  $ if $\mathrm{w}_{n-1}(C)=(n-1)\overline{n}$\\
\\
$v_{X}$ with $\mathrm{w}(X)=\widetilde{f}_{n-1}(\mathrm{w}(C))$ if
$\varphi_{n-1}(\mathrm{w}(C))=1$ and $\mathrm{w}_{n-1}(C)$ do not contain a
factor $\overline{n}n$\\
\\
$0$ otherwise
\end{tabular}
\right.  \label{fi_D}%
\end{equation}
where for any $p\in\mathbb{N}$, $(n\overline{n})^{p}$ and $(\overline{n}%
n)^{p}$ are respectively the words obtained by repeating $p$ times the factors
$n\overline{n}$ and $\overline{n}n$.

Given any column $C$ of type $D,$ we denote by $C_{n\longleftrightarrow
\overline{n}}$ the column obtained by permuting the letters $n$ and
$\overline{n}$ in $C.\;$The action of $f_{n}$ on $v_{C}$ may be obtained from
that of $f_{n-1}$ by using the involution of $W(\omega_{p}^{D})$ defined by%
\[
\Phi:v_{C}\longmapsto v_{C_{n\longleftrightarrow\overline{n}}}.
\]
We have%
\begin{equation}
f_{n}(v_{C})=\Phi f_{n-1}\Phi(v_{C}). \label{fn_D}%
\end{equation}

\noindent\textbf{Remark}: If we use the module $\mathcal{W}(\omega_{p})$
instead of $W(\omega_{p}),$ the coefficients of the decomposition of
$f_{i}(v_{C})$ on the basis $\{v_{C},$ $C\in\mathbf{Ca}(n,p)\}$ are not
necessary in $\mathbb{Z}[q]$.\ For example when $n=2$ one has by
(\ref{1_1bar_for_B}) and (\ref{1_1bar_for_D})%
\begin{equation}
f_{1}\left(  \
\begin{tabular}
[c]{|l|}\hline
$\mathtt{1}$\\\hline
$\mathtt{\bar{2}}$\\\hline
\end{tabular}
\ \right)  =%
\begin{tabular}
[c]{|l|}\hline
$\mathtt{2}$\\\hline
$\mathtt{\bar{2}}$\\\hline
\end{tabular}
+q^{2}%
\begin{tabular}
[c]{|l|}\hline
$\mathtt{1}$\\\hline
$\mathtt{\bar{1}}$\\\hline
\end{tabular}
=(1+q^{2})%
\begin{tabular}
[c]{|l|}\hline
$\mathtt{2}$\\\hline
$\mathtt{\bar{2}}$\\\hline
\end{tabular}
+\frac{q^{5}}{q^{4}-1}%
\begin{tabular}
[c]{|l|}\hline
$\mathtt{0}$\\\hline
$\mathtt{0}$\\\hline
\end{tabular}
\text{ for type }B \label{f1(12bar)B}%
\end{equation}
and
\begin{equation}
f_{1}\left(  \
\begin{tabular}
[c]{|l|}\hline
$\mathtt{1}$\\\hline
$\mathtt{\bar{2}}$\\\hline
\end{tabular}
\ \right)  =%
\begin{tabular}
[c]{|l|}\hline
$\mathtt{2}$\\\hline
$\mathtt{\bar{2}}$\\\hline
\end{tabular}
+q%
\begin{tabular}
[c]{|l|}\hline
$\mathtt{1}$\\\hline
$\mathtt{\bar{1}}$\\\hline
\end{tabular}
=\frac{1}{1-q^{2}}%
\begin{tabular}
[c]{|l|}\hline
$\mathtt{2}$\\\hline
$\mathtt{\bar{2}}$\\\hline
\end{tabular}
+\frac{q^{2}}{1-q^{2}}%
\begin{tabular}
[c]{|l|}\hline
$\mathtt{\bar{2}}$\\\hline
$\mathtt{2}$\\\hline
\end{tabular}
\text{ for type }D \label{f1(12bar)D}%
\end{equation}
where we have written for short $C$ in place of $v_{C}.$

\section{Computation of the canonical basis of $V(\lambda)\label{Sec_Algo}$}

\subsection{Global basis of $V(\omega_{p})$}

We can use Marsh's algorithm to compute the expansion of the global basis of
$V(\omega_{p})$ considered as a sub-module of $W(\omega_{p})$ on the basis
$\{v_{C},$ $C\in\mathbf{Ca}(n,p)\}.$

\begin{example}
Suppose $n=4$ and consider the admissible column $C=%
\begin{tabular}
[c]{|l|}\hline
$\mathtt{0}$\\\hline
$\mathtt{0}$\\\hline
$\mathtt{0}$\\\hline
$\mathtt{0}$\\\hline
\end{tabular}
$ of type $B$.\ Then $G(C)=f_{4}f_{3}f_{2}f_{1}f_{4}f_{3}f_{2}f_{4}f_{3}%
f_{4}v_{\omega_{3}}$ and we obtain:%
\[
G(C)=%
\begin{tabular}
[c]{|l|}\hline
$\mathtt{0}$\\\hline
$\mathtt{0}$\\\hline
$\mathtt{0}$\\\hline
$\mathtt{0}$\\\hline
\end{tabular}
+q(1+q^{6})%
\begin{tabular}
[c]{|l|}\hline
$\mathtt{4}$\\\hline
$\mathtt{0}$\\\hline
$\mathtt{0}$\\\hline
$\mathtt{\bar{4}}$\\\hline
\end{tabular}
+q(1-q^{4})%
\begin{tabular}
[c]{|l|}\hline
$\mathtt{3}$\\\hline
$\mathtt{0}$\\\hline
$\mathtt{0}$\\\hline
$\mathtt{\bar{3}}$\\\hline
\end{tabular}
+q(1+q^{2})%
\begin{tabular}
[c]{|l|}\hline
$\mathtt{2}$\\\hline
$\mathtt{0}$\\\hline
$\mathtt{0}$\\\hline
$\mathtt{\bar{2}}$\\\hline
\end{tabular}
+q^{2}(1+q^{2})(1-q^{4})%
\begin{tabular}
[c]{|l|}\hline
$\mathtt{3}$\\\hline
$\mathtt{4}$\\\hline
$\mathtt{\bar{4}}$\\\hline
$\mathtt{\bar{3}}$\\\hline
\end{tabular}
+q^{2}(1+q^{2})^{2}%
\begin{tabular}
[c]{|l|}\hline
$\mathtt{2}$\\\hline
$\mathtt{4}$\\\hline
$\mathtt{\bar{4}}$\\\hline
$\mathtt{\bar{2}}$\\\hline
\end{tabular}
\]
where we have written for short $C$ in place of $v_{C}$.
\end{example}

\noindent Note that the columns $X$ and the coefficients $d_{X}\in
\mathbb{Z}[q]$ such that $G(C)=\underset{X\in\mathbf{C}(n,p)}{\sum}%
d_{X}(q)v_{X}$ are in general more complex as in the symplectic case \cite{L1}
and it seens to be impossible to describe them in a simple combinatoric way.

\subsection{The representation $W(\lambda)$}

To make our notation homogeneous write for the spin representations
$W(\Lambda_{n}^{B})=V(\Lambda_{n}^{B})$, $W(\Lambda_{n}^{D})=V(\Lambda_{n}%
^{D})$ and $W(\Lambda_{n-1}^{D})=V(\Lambda_{n-1}^{D})$.

\noindent Consider $\lambda\in P^{+}$ a dominant weight and write
$\lambda=\Lambda+\lambda^{\prime}$ accordingly to (\ref{decom_lambda}). Since
$\lambda^{\prime}\in\Omega_{+},$ there exists $(\mu_{1},...,\mu_{n}%
)\in\mathbb{N}^{n}$ such that $\lambda^{\prime}=\underset{p=1}{\overset
{n}{\sum}}\mu_{p}\omega_{p}$.\ Then we set%
\begin{equation}
W(\lambda^{\prime})=W(\omega_{1})^{\otimes\mu_{1}}\otimes\cdot\cdot
\cdot\otimes W(\omega_{n})^{\otimes\mu_{n}}. \label{W(lambda)}%
\end{equation}
The natural basis of $W(\lambda^{\prime})$ consists of the tensor products
$v_{C_{r}}\otimes\cdot\cdot\cdot\otimes v_{C_{1}}$ of basis vectors $v_{C}$ of
the previous section appearing in (\ref{W(lambda)}).\ It is naturally indexed
by the tabloids of shape $Y(\lambda)$ by setting $v_{\tau}=v_{C_{r}}%
\otimes\cdot\cdot\cdot\otimes v_{C_{1}}$.

\noindent When $\lambda\neq\lambda^{\prime}$ we set
\[
W(\lambda)=W(\Lambda)\otimes W(\lambda^{\prime}).
\]
The natural basis of $W(\lambda)$ consists of the vectors $v_{\frak{C}}\otimes
v_{\tau}^{\prime}$ where $\tau^{\prime}$ is a tabloid of shape $Y(\lambda
^{\prime})$.\ It is naturally indexed by the tabloids of shape $Y(\lambda)$ by
setting $v_{\tau}=v_{\frak{C}}\otimes v_{\tau^{\prime}}.$

\bigskip

Let $L_{\lambda}$ be the $A$-submodule of $W(\lambda)$ generated by the
vectors $v_{\tau},$ $\tau\in\mathbf{T}(n,\lambda)$. We identify the image of
the vector $v_{\tau}$ by the projection $\pi_{\lambda}:L_{\lambda}\rightarrow
L_{\lambda}/qL_{\lambda}$ with the word \textrm{w(}$\tau)$. The pair
$(L_{\lambda},B_{\lambda}=\{\mathrm{w(}\tau)\left|  {}\right.  \tau
\in\mathbf{T}(n,\lambda)\})$ is then a crystal basis of $W(\lambda)$. Indeed
by Lemma \ref{lem_base_cryst_Wp}, it is the tensor product of the crystal
bases of the representations $W(\omega_{p})$ or $W(\Lambda)$ occurring in
$W(\lambda)$. Set%
\[
\left\{
\begin{tabular}
[c]{l}%
$v_{\lambda}=v_{\omega_{1}}^{\otimes\mu_{1}}\otimes\cdot\cdot\cdot\otimes
v_{\omega_{p}}^{\otimes\mu_{p}}\text{ if }\lambda\in\Omega^{+}$\\
$v_{\lambda}=v_{\Lambda}\otimes v_{\lambda^{\prime}}\text{ otherwise}$%
\end{tabular}
\right.  .
\]
We identify $V(\lambda)$ with the submodule of $W(\lambda)$ generated by
$v_{\lambda}$. Then, with the above notations, $v_{\lambda}=v_{T_{\lambda}}$
where $T_{\lambda}$ is the orthogonal tableau of shape $Y(\lambda)$ whose
$k$-th row is filled by letters $k$ for $k=1,...,n-1$ and $n$-th row by
letters $n$ except for $\lambda\in\Omega_{+}^{D}$ and $Y(\lambda)=(Y_{\lambda
},-)$ when it is filled by letters $\overline{n}$.\ By Theorem 4.2 of
\cite{Ka2}, we know that
\[
B(\lambda)=\{\widetilde{f}_{i_{1}}^{a_{1}}\cdot\cdot\cdot\widetilde{f}_{i_{r}%
}^{a_{r}}\mathrm{w}(T_{\lambda});\text{ }i_{1},...,i_{r}=1,...,n;\text{ }%
a_{1},...,a_{r}>0\}-\{0\}.
\]
The actions of $\widetilde{e}_{i}$ and $\widetilde{f}_{i},$ $i=1,...,n$ on
each $\mathrm{w(}\tau)\in B_{\lambda}$ are identical to those obtained by
considering $\mathrm{w(}\tau)$ as a vertex of $G_{n}$ or $\frak{G}_{n}$ since
it is true on each $B_{p}$. Hence, Definition \ref{defKN}, $B(\lambda
)=\{\mathrm{w}(T);$ $T\in\mathbf{OT}(n,\lambda)\}$. For each vector $G(T)$ of
the canonical basis of $V(\lambda)$ we will have:%
\[
G(T)\equiv v_{T}\operatorname{mod}qL_{\lambda}.
\]
The aim of this section is to describe an algorithm computing the
decomposition of the canonical basis $\{G(T);T\in\mathbf{OT}(n,\lambda)\}$
onto the basis $\{v_{\tau};\tau\in\mathbf{T}(n,\lambda)\}$ of $W(\lambda)$.

In the sequel we will need a total order on the readings of the tabloids. Let
$w_{1}=x_{1}\cdot\cdot\cdot x_{l}$ and $w_{2}=y_{1}\cdot\cdot\cdot y_{l}$ be
two distinct vertices of $G_{n}$ with the same length and $k$ the lowest
integer such that $x_{k}\neq y_{k}.$ When $w_{1}$ and $w_{2}$ are words of
$\mathcal{B}_{n}^{\ast},$ we write $w_{1}\trianglelefteq^{B}w_{2}$ if
$x_{k}\preceq y_{k}$ in $\mathcal{B}_{n}$ and $w_{1}\vartriangleright^{B}%
w_{2}$ otherwise that is, $\trianglelefteq^{B}$is the lexicographic order on
$\mathcal{B}_{n}^{\ast}$. When $w_{1}$ and $w_{2}$ are words of $\mathcal{D}%
_{n}^{\ast},$ we write $w_{1}\trianglelefteq^{D}w_{2}$ if $x_{k}\preceq y_{k}$
in $\mathcal{B}_{n}$ (indeed $\mathcal{D}_{n}$ is not totally ordered but the
letters of $\mathcal{D}_{n}$ are letters of $\mathcal{B}_{n})$.\ We can extend
the order $\trianglelefteq$ to vertices of $\frak{G}_{n}$ by setting
$w_{1}\frak{C}_{1}\trianglelefteq w_{2}\frak{C}_{2}\Longleftrightarrow
w_{1}\mathrm{w}(C_{\frak{C}_{1}})\trianglelefteq w_{1}\mathrm{w}%
(C_{\frak{C}_{2}}).$ We endow the set $\mathbf{T}(n,\lambda)$ with the total
order:%
\[
\tau_{1}\trianglelefteq\tau_{2}\Longleftrightarrow\mathrm{w(}\tau
_{1})\trianglelefteq\mathrm{w(}\tau_{2}).
\]
In fact $\trianglelefteq$ is a lexicographic order defined on the readings of
the tabloids of $\mathbf{T}(n,\lambda)$ such that%
\[
\left\{
\begin{tabular}
[c]{c}%
$x\trianglelefteq\widetilde{f}_{i}(x)$ for any letter $x$ with $\varphi
_{i}(x)\neq0$\\
$\frak{C}\trianglelefteq\widetilde{f}_{i}(\frak{C})$ for any spin column
$\frak{C}$ with $\varphi_{i}(\frak{C})\neq0$%
\end{tabular}
\right.  .
\]

We are going to compute the canonical basis $\{G(T);T\in\mathbf{OT}%
(n,\lambda)\}$ in two steps. First we obtain an intermediate basis
$\{A(T);T\in\mathbf{OT}(n,\lambda)\}$ which is fixed by the involution
$\overline{\text{%
\begin{tabular}
[c]{l}%
\ \
\end{tabular}
}}$ (condition (\ref{cond_invo})). When $T=C$ is an admissible column, $A(C)$
is the output of Marsh's algorithm.\ This implies that $A(C)=G(C)$.\ In the
general case, we have to correct $A(T)$ in order to verify condition
(\ref{cond_cong}). This second step is easy because we can prove that the
transition matrix from $\{A(T);T\in\mathbf{OT}(n,\lambda)\}$ to $\{G(T);T\in
\mathbf{OT}(n,\lambda)\}$ is unitriangular once the orthogonal tableaux and
the tabloids are ordered by $\trianglelefteq$.

We start with a general Lemma analogous to Lemma 4.1 of \cite{L1}.

\begin{lemma}
\label{Lem_F()}Let $v\in V(\lambda)$ be a vector of the type
\begin{equation}
v=f_{i_{1}}^{(r_{1})}\cdot\cdot\cdot f_{i_{s}}^{(r_{s})}v_{\lambda}
\label{Type_F()}%
\end{equation}
where $(i_{1},...,i_{s})$ and $(r_{1},...,r_{s})$ are two sequences of
integers.\ Then the coordinates of $v$ on the basis $\{v_{\tau};\tau
\in\mathbf{T}(n,\lambda)\}$ belong to $\mathbb{Z}[q,q^{-1}]$.
\end{lemma}

\begin{proof}
It suffices to prove that for any tabloid $\tau=C_{1}\cdot\cdot\cdot C_{s},$
any $i=1,...,n$ and any integer $m,$ the coordinates of $f_{i}^{(m)}v_{\tau}$
on the basis $\{v_{\tau};\tau\in\mathbf{T}(n,\lambda)\}$ belong to
$\mathbb{Z}[q,q^{-1}].$ We proceed by induction on $s$. If $s=1,$ the result
follows from the description of the action of the Chevalley generators on the
vectors $v_{C}$ given in \ref{sub_sec_formulas}. Set $\tau^{\prime}=C_{1}%
\cdot\cdot\cdot C_{s-1}$.\ Then $v_{\tau}=v_{C_{s}}\otimes v_{\tau^{\prime}}.$
By Lemma \ref{lem_sl2}, we have%
\[
f_{i}^{(m)}(v_{C_{s}}\otimes v_{\tau^{\prime}})=\overset{m}{\underset
{k=0}{\sum}}q_{i}^{(m-k)(a-k)}\ f^{(k)}(v_{C_{s}})\otimes f^{(m-k)}%
(v_{\tau^{\prime}})
\]
where $t_{i}(v_{C_{s}})=q_{i}^{a}v_{C_{s}}.$ By induction the coordinates of
$f^{(k)}(v_{C_{s}})$ and $f^{(m-k)}(v_{\tau^{\prime}})$ respectively on the
bases $\{v_{C};C\in\mathbf{C}(n,p)\}$ and $\{v_{\tau};\tau\in\mathbf{T}%
(n,\lambda)\}$ belong to $\mathbb{Z}[q,q^{-1}].\;$So it is true for
$f_{i}^{(m)}(v_{\tau}).$
\end{proof}

Note that this lemma is false if we use $\mathcal{W}(\omega_{p})$ instead of
$W(\omega_{p})$ since $\underset{C\in\mathbf{Ca}(n,p)}{\sum}\mathbb{Z}%
[q,q^{-1}]v_{C}$ is not stabilized by the $f_{i}$'s.

\subsection{The basis $A(T)$}

The basis $\{A(T)\}$ will be a monomial basis, that is, a basis of the form
\begin{equation}
A(T)=f_{i_{1}}^{(r_{1})}\cdot\cdot\cdot f_{i_{m}}^{(r_{m})}v_{\lambda}.
\label{A(T)_monom}%
\end{equation}
By Lemma \ref{Lem_F()}, the coordinates of $A(T)$ on the basis $\{v_{\tau}\}$
of $W(\lambda)$ belongs to $\mathbb{Z[}q,q^{-1}]$. To find the two sequences
of integers $(i_{1},...,i_{m})$ and $(r_{1},...,r_{m})$ associated to $T$, we
proceed as follows.

\noindent Suppose first $\lambda\in\Omega_{+}$ and consider $T=C_{1}\cdot
\cdot\cdot C_{s}\neq T_{\lambda}\in\mathbf{OT}(n,\lambda)$.\ Let $C_{k}$ be
the rightmost column of $T$ such that $\mathrm{w(}C_{k})$ is not a highest
weight vertex.\ We define $i_{1}\in\{1,...,n\}$ from $C_{k}$ as we have done
in \ref{sub-sec_M_algo}. If $k=1$ we set $l=1.\;$If $k>1$ and $f_{i_{1}%
}v_{C_{k}}\not =0$ or $\widetilde{e}_{i_{1}}(\mathrm{w}(C_{k-1}))=0$ we set
$l=k.\;$Otherwise let $l$ be the lowest integer $l<k$ satisfying the two
conditions%
\begin{equation}
\left\{
\begin{tabular}
[c]{l}%
$\mathrm{(i):}$ $f_{i_{1}}v_{C_{j}}=0$ for $j=l+1,...,k$\\
$\mathrm{(ii):}$ $\widetilde{e}_{i_{1}}(\mathrm{w}(C_{j}))\neq0$ for
$j=l,...,k$%
\end{tabular}
\right.  . \label{cond}%
\end{equation}
Set $\varepsilon_{i_{1},j}=\varepsilon_{i_{1}}((\mathrm{w}(C_{j}))$ for
$j=l,...,k$ and $r_{1}=\overset{k}{\underset{j=l}{\sum}}\varepsilon_{i_{1}j}%
$.\ Write $T_{1}$ for the tabloid obtained by changing in $T$ each column
$C_{j},$ $j=l,...,k$ into the column of reading $\widetilde{e}_{i_{1}%
}^{\varepsilon_{i_{1},j}}(\mathrm{w}(C_{j})).$

\noindent Now suppose $\lambda\notin\Omega_{+}$ and consider $T=\frak{C}%
C_{1}\cdot\cdot\cdot C_{s}\neq T_{\lambda}\in\mathbf{OT}(n,\lambda).\;$Write
$T^{\prime}=C_{1}\cdot\cdot\cdot C_{s}$.\ If $T^{\prime}$ is of highest
weight, we have $G(T)=v_{T}.$ Indeed the vectors $v_{\tau}$ appearing in the
decomposition of $G(T)$ on the basis $\{v_{\tau},$ $\tau\in\mathbf{T}%
(n,\lambda)\}$ have the same weight $\mu$ than $T$ and $T$ is then the unique
tabloid of shape $\lambda$ and weight $\mu.\;$So we can suppose that
$T^{\prime}$ is not of highest weight and define $k,$ $l,$ $i_{1},$
$r_{1}^{\prime}$ and $T_{1}^{\prime}$ from $T^{\prime}$ as we have done in the
case $\lambda\in\Omega^{+}$.\ Finally we set%
\begin{align*}
T_{1}  &  =(\widetilde{e}_{i_{1}}\frak{C})T_{1}^{\prime}\text{ and }%
r_{1}=r_{1}^{\prime}+1\text{ if }l=1\text{, }f_{i_{1}}(v_{C_{l}})=0\text{ and
}\varepsilon_{i_{1}}(\frak{C})=1,\\
T_{1}  &  =\frak{C}T_{1}^{\prime}\text{ and }r_{1}=r_{1}^{\prime}\text{
otherwise.}%
\end{align*}

\begin{lemma}
\label{T1_Tab}$T_{1}\in\mathbf{OT}(n,\lambda).$
\end{lemma}

\begin{proof}
Suppose first that $\lambda\in\Omega_{+}.$ Set $T_{1}=D_{1}\cdot\cdot\cdot
D_{s}$.\ Then by construction of $T_{1},$ we have $D_{i}=C_{i}$ for
$i\notin\{l,...,k\}.\;$The columns $C_{i}$ with $i=k+1,...,s$ are of highest
weight.\ One has $\mathrm{w}(D_{k}\cdot\cdot\cdot D_{s})=\widetilde{e}_{i_{1}%
}^{\varepsilon_{i_{1},k}}(\mathrm{w(}C_{k}\cdot\cdot\cdot C_{s}))$ by
(\ref{TENS2}) because $\varepsilon_{i_{1},k}=\varepsilon_{i_{1}}%
\mathrm{w}(C_{k}\cdot\cdot\cdot C_{s})$ and $\varphi_{i_{1}}\mathrm{w(}%
C_{j})=0$ for $j=k+1,...,s$ (otherwise we would have $\widetilde{e}_{i_{1}%
}(\mathrm{w(}C_{k}))=0$ since the letters $1,...,i_{1}$ would occur in $C_{k}%
$).\ So $D_{k}\cdot\cdot\cdot D_{s}$ is an orthogonal tableau.

\noindent If $l=1,$ $\mathrm{w}(D_{1}\cdot\cdot\cdot D_{k})=\widetilde
{e}_{i_{1}}^{r_{1}}(\mathrm{w}(C_{1}\cdot\cdot\cdot C_{k}))$ since
$\varphi_{i_{1}}(\mathrm{w(}C_{j}))=0$ for $j=2,...,k.\;$So $D_{1}\cdot
\cdot\cdot D_{k}$ is an orthogonal tableau.$\;$The lemma is true since
$D_{k}\cdot\cdot\cdot D_{s}$ is also an orthogonal tableau. If $l>1,$ we have
either $f_{i_{1}}v_{C_{l}}\neq0$ and $\widetilde{e}_{i_{1}}(\mathrm{w}%
(C_{l}))\neq0$ or $\widetilde{e}_{i_{1}}(\mathrm{w}(C_{l-1}))=0.$

\noindent When $f_{i_{1}}v_{C_{l}}\neq0$, by considering the configurations
given in \ref{sub_sec_formulas} for $f_{i_{1}}v_{C_{l}}\neq0$ and
$\widetilde{e}_{i_{1}}(\mathrm{w}(C_{l}))\neq0$, we see that $\varepsilon
_{i_{1},l}=1$ and $\varphi_{i_{1}}\mathrm{w(}D_{l})=2$.\ More precisely one
has:%
\begin{align}
\mathrm{w}_{i_{1}}\mathrm{(}D_{l})  &  =\left\{
\begin{tabular}
[c]{l}%
$i_{1}\overline{i_{1}+1}\text{ if }i_{1}=1,...,n-1$\\
$n0^{r}$ if $i_{1}=n$%
\end{tabular}
\right.  \text{ when }D_{l}\text{ is of type }B\label{HP1}\\
\mathrm{w}_{i_{1}}\mathrm{(}D_{l})  &  =\left\{
\begin{tabular}
[c]{l}%
$i_{1}\overline{i_{1}+1}\text{ if }i_{1}=1,...,n-2$\\
$(n-1)(\overline{n}n)^{r}\overline{n}$ if $i_{1}=n-1$\\
$(n-1)(n\overline{n})^{r}n$ if $i_{1}=n$%
\end{tabular}
\right.  \text{ when }D_{l}\text{ is of type }D \label{HP2}%
\end{align}
In all cases, a simple computation from Definition 3.1.5 of \cite{L2(RS)}
shows that $lC_{l}=lD_{l}$ (see Remark following Definition \ref{defKN}%
).\ Hence $rC_{l-1}lD_{l}$ is an orthogonal tableau.

\noindent When $\widetilde{e}_{i_{1}}(\mathrm{w}(C_{l-1}))=0,$ we have
$\widetilde{e}_{i_{1}}^{\varepsilon_{i_{i},l}}(\mathrm{w}(C_{l})\mathrm{w}%
(C_{l-1}))=\mathrm{w}(D_{l})\mathrm{w}(C_{l-1})$ by (\ref{TENS2}).\ So
$C_{l-1}D_{l}$ is an orthogonal tableau.

\noindent Finally we obtain that $D_{1}\cdot\cdot\cdot D_{l}=C_{1}\cdot
\cdot\cdot C_{l-1}D_{l}$ is an orthogonal tableau.$\;$Moreover $\mathrm{w(}%
D_{l}\cdot\cdot\cdot D_{k})=\widetilde{e}_{i_{1}}^{r_{1}}\mathrm{w(}%
(C_{l}\cdot\cdot\cdot C_{k}))$ because $\varphi_{i_{1}}\mathrm{w(}C_{j})=0$
for $j=l+1,...,k$ so is the reading of an orthogonal tableau. Hence $T_{1}$ is
an orthogonal tableau since $D_{k}\cdot\cdot\cdot D_{s}$ is an orthogonal tableau.

Now suppose that $\lambda\notin\Omega_{+}.$ Set $T_{1}=\frak{D}D_{1}\cdot
\cdot\cdot D_{s}.$ We know by the arguments above that $T_{1}^{\prime}%
=D_{1}\cdot\cdot\cdot D_{s}$ is an orthogonal tableau. So it suffices to prove
that $\frak{D}D_{1}$ is an orthogonal tableau. If $l>1,$ $\frak{D=C}$ and
$D_{1}=C_{1}$. Now suppose $l=1.\;$If $f_{i_{1}}(v_{C_{1}})=0$ and
$\varepsilon_{i_{1}}(\frak{C})=1$ then $\mathrm{w(}\frak{D}D_{1}%
)=\widetilde{e}_{i_{1}}^{\varepsilon_{i_{1},l}+1}\mathrm{w(}\frak{C}C_{1})$ is
the reading of an orthogonal tableau.\ If $f_{i_{1}}(v_{C_{1}})=0$ and
$\varepsilon_{i_{1}}(\frak{C})=0$ then $\mathrm{w(}\frak{D}D_{1}%
)=\widetilde{e}_{i_{1}}^{\varepsilon_{i_{1},1}}\mathrm{w(}\frak{C}C_{1})$ is
the reading of an orthogonal tableau. Finally if $f_{i_{1}}(v_{C_{1}})\neq0$
we have seen that $lC_{1}=lD_{1}.\;$Hence $\mathrm{w(}\frak{D}D_{1}%
)=\mathrm{w(}\frak{C}D_{1})$ is again the reading of an orthogonal tableau.
\end{proof}

\begin{remark}
\label{remark}For $j=l,...,k,$ $\widetilde{e}_{i_{1}}(\mathrm{w}(D_{j}))=0$
and $\widetilde{f}_{i_{1}}(\mathrm{w}(D_{j}))\neq0$.\ By considering the
configurations that can occur in $D_{j},$ this implies that $e_{i_{1}%
}(v_{D_{j}})=0$ and $f_{i_{1}}^{(\varepsilon_{i_{1},j})}v_{D_{j}}=v_{C_{j}}$
for $j=l+1,...,k.$
\end{remark}

Once the tableau $T_{1}$ defined, we do the same with $T_{1}$ getting a new
orthogonal tableau $T_{2}$ and a new integer $i_{2}$. And so on until the
tableau $T_{s}$ obtained is equal to $T_{\lambda}$. Notice that we can not
write \textrm{w(}$T_{1})=\widetilde{e}_{i_{1}}^{r_{1}}\mathrm{w(}T)$ in
general, that is, our algorithm does not provide a path in the crystal graph
$B(\lambda)$ joining the vertex $\mathrm{w(}T)$ to the vertex of highest
weight $\mathrm{w(}T_{\lambda})$.

\begin{example}
For $T=%
\begin{tabular}
[c]{|l|ll}\hline
$\mathtt{2}$ & $\mathtt{2}$ & \multicolumn{1}{|l|}{$\mathtt{3}$}\\\hline
$\mathtt{0}$ & $\mathtt{\bar{3}}$ & \multicolumn{1}{|l}{}\\\cline{1-1}%
\cline{1-2}%
$\mathtt{0}$ &  & \\\cline{1-1}%
\end{tabular}
$ of type $B$ and $n=3$, we obtain successively\vspace{0.5cm}%

\begin{tabular}
[c]{|l|ll}\hline
$\mathtt{2}$ & $\mathtt{2}$ & \multicolumn{1}{|l|}{$\mathtt{2}$}\\\hline
$\mathtt{0}$ & $\mathtt{\bar{3}}$ & \multicolumn{1}{|l}{}\\\cline{1-1}%
\cline{1-2}%
$\mathtt{0}$ &  & \\\cline{1-1}%
\end{tabular}
,
\begin{tabular}
[c]{|l|ll}\hline
$\mathtt{1}$ & $\mathtt{1}$ & \multicolumn{1}{|l|}{$\mathtt{1}$}\\\hline
$\mathtt{0}$ & $\mathtt{\bar{3}}$ & \multicolumn{1}{|l}{}\\\cline{1-1}%
\cline{1-2}%
$\mathtt{0}$ &  & \\\cline{1-1}%
\end{tabular}
,
\begin{tabular}
[c]{|l|ll}\hline
$\mathtt{1}$ & $\mathtt{1}$ & \multicolumn{1}{|l|}{$\mathtt{1}$}\\\hline
$\mathtt{3}$ & $\mathtt{3}$ & \multicolumn{1}{|l}{}\\\cline{1-1}\cline{1-2}%
$\mathtt{0}$ &  & \\\cline{1-1}%
\end{tabular}
,
\begin{tabular}
[c]{|l|ll}\hline
$\mathtt{1}$ & $\mathtt{1}$ & \multicolumn{1}{|l|}{$\mathtt{1}$}\\\hline
$\mathtt{2}$ & $\mathtt{2}$ & \multicolumn{1}{|l}{}\\\cline{1-1}\cline{1-2}%
$\mathtt{0}$ &  & \\\cline{1-1}%
\end{tabular}
and
\begin{tabular}
[c]{|l|ll}\hline
$\mathtt{1}$ & $\mathtt{1}$ & \multicolumn{1}{|l|}{$\mathtt{1}$}\\\hline
$\mathtt{2}$ & $\mathtt{2}$ & \multicolumn{1}{|l}{}\\\cline{1-1}\cline{1-2}%
$\mathtt{3}$ &  & \\\cline{1-1}%
\end{tabular}
.%

\[
A(T)=f_{2}f_{1}^{(3)}f_{3}^{(3)}f_{2}^{(2)}f_{3}T_{\lambda}.
\]
\end{example}

\begin{proposition}
\label{prop_A(T)}The expansion of $A(T)$ on the basis $\{v_{\tau};\tau
\in\mathbf{T}(n,\lambda)\}$ of $W(\lambda)$ is of the form
\[
A(T)=\underset{\tau}{\sum}\alpha_{\tau,T}(q)v_{\tau}%
\]

where the coefficients $\alpha_{\tau,T}(q)$ satisfy:

\textrm{(i)}: $\alpha_{\tau,T}(q)\neq0$ only if $\tau$ and $T$ have the same weight,

\textrm{(ii)}: $\alpha_{\tau,T}(q)\in\mathbb{Z}[q,q^{-1}]$ and $\alpha_{T,T}(q)=1,$

\textrm{(iii)}: $\alpha_{\tau,T}(q)\neq0$ only if $\tau\trianglelefteq T.$
\end{proposition}

\begin{proof}
(of Proposition \ref{prop_A(T)})

We only sketch the proof for $\lambda\in\Omega_{+}.\;$When $\lambda
\notin\Omega_{+}$ the arguments are essentially the same up to minor
modifications due to the spin column contained in the tabloids of shape
$\lambda.$ Consider a tabloid $\tau\in\mathbf{T}(n,\lambda)$ with $\lambda
\in\Omega_{+}$ and $i\in\{1,...n\}$.\ We call ``$i$-substitution in $\tau$''
the substitution of any letter $a\in\tau$ such that $\varphi_{i}(a)\neq0$ by
the letter $\widetilde{f}_{i}(a)$ preserving the structure of tabloid.

\textrm{(i) }is a straightforward consequence of the definition of $A(T)$. By
Lemma \ref{Lem_F()}, we know that $\alpha_{\tau,T}(q)\in\mathbb{Z}[q,q^{-1}].$
The proposition will be proved by induction if we show that \textrm{(ii)
}and\textrm{ (iii) }hold for $T$ as soon as they hold for $T_{1}$ with
$A(T)=f_{i_{1}}^{(r_{1})}A(T_{1})$ (the notations are those of Lemma
\ref{T1_Tab}). If the vector $v_{\tau}$ occurs in $A(T)$, the tabloid $\tau$
is obtained from a tabloid $\tau_{1}$ labelling a vector $v_{\tau_{1}}$
occurring in $A(T_{1})$ after $r_{1}$ $i_{1}$-substitutions in $\tau_{1}$.

When $\tau_{1}=T_{1}$ and $\tau=T,$ the definition of $T_{1}$ implies that
$\mathrm{w}(T)$ is obtained from $\mathrm{w}(T_{1})$ by executing the $r_{1}$
leftmost possible $i_{1}$-substitutions in $\mathrm{w}(T_{1}).$ Hence $v_{T}$
appears in $f_{i_{1}}^{(r_{1})}v_{T_{1}}$ with a non zero coefficient. Now
suppose that there exists $\tau_{1}\neq T_{1}$ such that $v_{\tau_{1}}$
appears in $A(T_{1})$ and $v_{T}$ appears in $f_{i_{1}}^{(r_{1})}v_{\tau_{1}}%
$. Write $\mathrm{w}(T)=z_{1}\cdot\cdot\cdot z_{r},$ $\mathrm{w}(T_{1}%
)=y_{1}\cdot\cdot\cdot y_{r}$ and $\mathrm{w}(\tau_{1})=x_{1}\cdot\cdot\cdot
x_{r}.$ Let $w$ be the factor of the words \textrm{w}$(\tau_{1})$ and
\textrm{w}$(T_{1})$ of maximal length such that there exist two words
$u,u^{\prime}$ and two letters $x_{q}\neq y_{q}$ satisfying:%
\begin{equation}
\mathrm{w}(\tau_{1})=wx_{q}u\text{ and }\mathrm{w}(T_{1})=wy_{q}u^{\prime}.
\label{dec}%
\end{equation}
We must have $x_{q}\vartriangleleft y_{q}$ because $\tau_{1}\vartriangleleft
T_{1}$. The letter $x_{q}$ is necessarily modified when $T$ is obtained from
$\tau_{1}$. Otherwise we have $x_{q}\in$ \textrm{w}$(T)$. But \textrm{w}$(T)$
can also be computed from \textrm{w}$(T_{1}).$ So the letter of \textrm{w}%
$(T)$ occurring at the same place than the letter $y_{q}$ in \textrm{w}%
$(T_{1})$ is $\trianglerighteq y_{q}$: it can not be $x_{q}$. This implies
that $x_{q}=\widetilde{e}_{i_{1}}(y_{q}).$

\noindent Write $\tau_{1}=E_{1}\cdot\cdot\cdot E_{s}$ where $E_{j}$ is a
column for $j=1,...,s$.\ Suppose that $x_{q}$ appear in the $p$-th column
$E_{p}$ of $\tau_{1}.$ Then $y_{q}\in D_{p}$ and $z_{q}=y_{q}\in C_{p}$ (the
notations are those of Lemma \ref{T1_Tab}, that is $T=C_{1}\cdot\cdot\cdot
C_{s}$ and $T_{1}=D_{1}\cdot\cdot\cdot D_{s}$).

\noindent Suppose $C_{p}\neq D_{p}$.$\;$Then $\mathrm{w}(D_{p})$ is obtained
by applying $\widetilde{e}_{i_{1}}^{\varepsilon}$ with $\varepsilon
=\varepsilon_{i_{1}}(\mathrm{w}(C_{p}))$ to $\mathrm{w}(C_{p})$.\ Hence
$C_{p}$ contains at least a letter $z_{h}$ which is changed into
$y_{h}=\widetilde{e}_{i_{1}}(z_{h})$ to obtain $D_{p}.\;$Necessarily $h\neq q$
because $z_{q}=y_{q}\in D_{p}$.\ So $C_{p}$ contain the two letters
$z_{q}=y_{q}$ and $z_{h}$.\ Then the diagram obtained by changing in $C_{p}$
$z_{q}$ and $z_{h}$ into $\widetilde{e}_{i_{1}}(z_{q})=x_{q}$ and
$\widetilde{e}_{i_{1}}(z_{h})=y_{h}$ is the column of reading $\widetilde
{e}_{i_{1}}^{2}\mathrm{w}(C_{p})$.\ Hence $\varepsilon=2$ and $\mathrm{w}%
(D_{p})=\widetilde{e}_{i_{1}}^{2}\mathrm{w}(C_{p})$ by definition of
$D_{p}.\;$It means that $y_{q}\notin D_{p}$ and we obtain a contradiction.\ We
have proved that $C_{p}=D_{p}.$

\noindent The column $C_{p}$ is not of highest weight because it contains the
movable letter $y_{q}.\;$So $p\leq k$ (see (\ref{cond})).$\;$This implies that
$p<l$ because $C_{p}$ is not modified when we obtain $T_{1}$ from $T.\;$Hence
\textrm{w}$(T)$ is obtained from \textrm{w}$(T_{1})$ after $r_{1}$ $i_{1}%
$-substitutions occurring in $w$.\ We derive a contradiction because in this
case there are $r_{1}+1$ $i_{1}$-substitutions when $T$ is obtained from
$\tau_{1}$. We have proved that $v_{T}$ can only appear in $f_{i_{1}}%
^{(r_{1})}v_{T_{1}}$.

\noindent The vector $v_{T}$ can only appear when $f_{i_{1}}^{(r_{1})}$ is
applied to $v_{D_{k}}\otimes\cdot\cdot\cdot\otimes v_{D_{l}}$ in $v_{T_{1}%
}.\;$We have $e_{i_{1}}(v_{D_{j}})=0$ (see Remark \ref{remark}) and $f_{i_{1}%
}(v_{C_{j}})=0$ for $j=l+1,...,k.$ Set $\alpha=\varphi_{i_{1}}(\mathrm{w}%
(D_{l+1}\cdot\cdot\cdot D_{k}))$, $v_{D}=v_{D_{k}}\otimes\cdot\cdot
\cdot\otimes v_{D_{l+1}}$ and $v_{C}=v_{C_{k}}\otimes\cdot\cdot\cdot\otimes
v_{C_{l+1}}.$ Lemma \ref{lem_sl2} and an immediate induction shows that
$f_{i_{1}}^{(\alpha)}v_{D}=v_{C}.$ By using one more time this lemma we have%

\[
f^{(r_{1})}(v_{D}\otimes v_{D_{l}})=\overset{r_{1}}{\underset{m=0}{\sum}%
}q^{(r_{1}-m)(\alpha-m)}\ f^{(m)}(v_{D})\otimes f^{(r_{1}-m)}(v_{D_{l}}).
\]
In this sum the vector $v_{C}\otimes v_{C_{l}}$ only appear for $m=\alpha$
because $v_{D}$ must be changed into $v_{C}$.\ Hence the coefficient of
$v_{C}\otimes f^{(r_{1}-\alpha)}(v_{D_{l}})=v_{C}\otimes v_{C_{l}}$ in the
above sum is $1$.\ So the coefficient of $v_{T}$ in $f_{i_{1}}^{(r_{1}%
)}v_{T_{1}}$ is equal to $1$ which proves $\mathrm{(ii)}$ since the
coefficient of $v_{T_{1}}$ in $A(T_{1})$ is $1$.

Consider $v_{\tau}$ appearing in $A(T)$ and suppose that the tabloid $\tau$ is
obtained from the tabloid $\tau_{1}$ such that $v_{\tau_{1}}$ appears in
$A(T_{1})$ after $r_{1}$ $i_{1}$-substitutions in $\tau_{1}$. Let
$\tau^{\prime}$ be the tabloid obtained by executing in \textrm{w}$(\tau_{1})$
the $r_{1}$ leftmost possible $i_{1}$-substitutions. We are going to prove
that $\tau^{\prime}\trianglelefteq T$ which implies the proposition because
$\tau\trianglelefteq\tau^{\prime}$. If $\tau_{1}=T_{1}$ then $\tau^{\prime}%
=T$. So we can suppose $\tau_{1}\neq T_{1}$ and decompose the words
$\mathrm{w}(\tau_{1})$, $\mathrm{w}(T_{1})$ as in (\ref{dec}) with
$x_{q}\vartriangleleft y_{q}$. If $\tau^{\prime}\trianglerighteq T,$ there is
a $i_{1}$-substitution executed in $w$ when $\tau^{\prime}$ is obtained from
$\tau_{1}$ which is not executed when $T$ is obtained from $T_{1}$. Denote by
$a$ the letter concerned by this substitution.\ We can write $w=w_{1}aw_{2}$
where $w_{1}$ and $w_{2}$ are words and%
\[
\mathrm{w}(\tau_{1})=w_{1}aw_{2}x_{q}u\text{ and }\mathrm{w}(T_{1}%
)=w_{1}aw_{2}y_{q}u^{\prime}.
\]
Then when we compute $T$ from $T_{1},$ the $r_{1}$ $i_{1}$-substitutions occur
in $w_{1}$ since $\mathrm{w}(T)$ is obtained from $\mathrm{w}(T_{1})$ by
executing the $r_{1}$ leftmost possible $i_{1}$-substitutions in
$\mathrm{w}(T_{1}).$ This contradicts the definition of $\tau^{\prime}$. So
\textrm{(iii) }is true.
\end{proof}

It follows from \textrm{(iii)} that the vectors $A(T)$ are linearly
independent in $V(\lambda)$. This implies that $\{A(T);T\in\mathbf{OT}%
(n,\lambda)\}$ is a $\mathbb{Q[}q]$-basis of $V(\lambda)$. Indeed by
Definition \ref{defKN}, $\dim V(\lambda)=\mathrm{card}(\mathbf{OT}%
(n,\lambda))$. As a consequence of (\ref{A(T)_monom}), we obtain
$\overline{A(T)}=A(T)$. Note that, by definition of Marsh's algorithm, the
bases $\{A(T)\}$ and $\{G(T)\}$ coincide when $\lambda=\omega_{i},$ $i=1,...,n.$

\subsection{From $A(T)$ to $G(T)$}

Let us write%
\begin{align*}
G(T)  &  =\underset{\tau\in\mathbf{T}^{B}(n,\lambda)}{\sum}d_{\tau,T}%
^{B}(q)\,v_{\tau}\text{ when }T\text{ is of type }B,\\
G(T)  &  =\underset{\tau\in\mathbf{T}^{D}(n,\lambda)}{\sum}d_{\tau,T}%
^{D}(q)\,v_{\tau}\text{ when }T\text{ is of type }D.
\end{align*}
We are going to describe a simple algorithm for computing the rectangular
matrix of coefficients
\[
D=[d_{\tau,T}(q)],\text{ \ \ \ \ }\tau\in\mathbf{T}(n,\lambda)\text{,
\ \ \ \ }T\in\mathbf{OT}(n,\lambda).
\]

\begin{lemma}
\label{G_on_tabloids}The coefficients $d_{\tau,T}(q)$ belong to $\mathbb{Q}%
[q]$. Moreover $d_{\tau,T}(0)=0$ if $\tau\neq T$ and $d_{T,T}=1$.
\end{lemma}

\begin{proof}
Recall that $\{G(T)\}$ is a basis of $V_{\mathbb{Q}}(\lambda)=U_{\mathbb{Q}%
}^{-}v_{\lambda}$. This implies that the vectors of this basis are
$\mathbb{Q[}q,q^{-1}]$-linear combinations of vectors of the type considered
in Lemma \ref{Lem_F()}.\ In particular $d_{\tau,T}\in\mathbb{Q}[q,q^{-1}%
]$.\ By condition (\ref{cond_cong}), $d_{\tau,T}(q)$ must be regular at $q=0$
hence $d_{\tau,T}(q)$ belongs to $\mathbb{Q[}q]$. Moreover
\[
d_{\tau,T}(q)\equiv\left\{
\begin{tabular}
[c]{l}%
$0$ $\operatorname{mod}q$ \ \ if $\tau\neq T$\\
$1$ $\operatorname{mod}q$ \ \ otherwise
\end{tabular}
\right.  .
\]
So the Lemma is true.
\end{proof}

Let us write
\begin{align}
G(T)  &  =\underset{S\in OT^{B}(n,\lambda)}{\sum}\beta_{S,T}^{B}%
(q)\,A(S)\text{ when }T\text{ is of type }B\label{G(T)_on_A(T)}\\
G(T)  &  =\underset{S\in OT^{D}(n,\lambda)}{\sum}\beta_{S,T}^{D}%
(q)\,A(S)\text{ when }T\text{ is of type }D\nonumber
\end{align}
the expansion of the basis $\{G(T)\}$ on the basis $\{A(T)\}$. We have the
following lemma analogous to Lemma 4.3 of \cite{L-T}:

\begin{lemma}
The coefficients $\beta_{S,T}(q)$ of (\ref{G(T)_on_A(T)}) satisfy:

$\mathrm{(i)}$: $\beta_{S,T}(q)=\beta_{S,T}(q^{-1}),$

$\mathrm{(ii)}$: $\beta_{S,T}(q)=0$ unless $S\trianglelefteq T,$

$\mathrm{(iii)}$: $\beta_{T,T}(q)=1$.
\end{lemma}

\begin{proof}
See proof of Lemma 4.3 in \cite{L-T}.
\end{proof}

Let $T_{\lambda}=T^{(1)}\vartriangleleft T^{(2)}\vartriangleleft\cdot
\cdot\cdot\vartriangleleft T^{(t)}$ be the sequence of tableaux of
$\mathbf{OT}(n,\lambda)$ ordered in increasing order. We have $G(T_{\lambda
})=A(T_{\lambda})$, i.e. $G(T^{(1)})=A(T^{(1)})$. By the previous lemma, the
transition matrix $M$ from $\{A(T)\}$ to $\{G(T)\}$ is upper unitriangular
once the two bases are ordered with $\trianglelefteq$. Since $\{G(T)\}$ is a
$\mathbb{Q}[q,q^{-1}]$ basis of $V_{\mathbb{Q}}(\lambda)$ and $A(T)\in
V_{\mathbb{Q}}(\lambda),$ the entries of $M$ are in $\mathbb{Q}[q,q^{-1}].$
Suppose by induction that we have computed the expansion on the basis
$\{v_{\tau};\tau\in\mathbf{T}(n,\lambda)\}$ of the vectors%
\[
G(T^{(1)}),...,G(T^{(i)})
\]
and that this expansion verifies $d_{\tau,T^{(p)}}(q)=0$ if $\tau
\vartriangleright T^{(p)}$ for $p=1,...,i$. The inverse matrix $M^{-1}$ is
also upper unitriangular with entries in $\mathbb{Q}[q,q^{-1}]$. So we can
write:%
\begin{equation}
G(T^{(i+1)})=A(T^{(i+1)})-\gamma_{i}(q)G(T^{(i)})-\cdot\cdot\cdot-\gamma
_{1}(q)G(T^{(1)})\text{.} \label{G_on_A}%
\end{equation}
It follows from condition (\ref{cond_invo}) and Proposition \ref{prop_A(T)}
that $\gamma_{m}(q)=\gamma_{m}(q^{-1})$ for $m=1,...,i$. By Lemma
\ref{G_on_tabloids}, the coordinate $d_{T^{(i)},T^{(i+1)}}(q)$ of
$G(T^{(i+1)})$ on the vector $v_{T^{(i)}}$ belongs to $\mathbb{Q}[q]$,
$d_{T^{(i)},T^{(i+1)}}(0)=0$ and the coordinate $d_{T^{(i)},T^{(i)}}(q)$ of
$G(T^{(i)})$ on the vector $v_{T^{(i)}}$ is equal to $1$. Moreover
$v_{T^{(i)}}$ can only occur in $A(T^{(i+1)})-\gamma_{i}(q)G(T^{(i)}).$ If
\[
\alpha_{T^{(i)},T^{(i+1)}}(q)=\underset{j=-r}{\overset{s}{\sum}}a_{j}q^{j}%
\in\mathbb{Z}[q,q^{-1}]
\]
then we will have
\[
\gamma_{i}(q)=\overset{0}{\underset{j=-r}{\sum}}a_{j}q^{j}+\underset
{j=1}{\overset{r}{\sum}}a_{-j}q^{j}\in\mathbb{Z}[q,q^{-1}].
\]
Next if the coefficient of $v_{T^{(i-1)}}$ in $A(T^{(i+1)})-\gamma
_{i}(q)G(T^{(i)})$ is equal to
\[
\underset{j=-l}{\overset{k}{\sum}}b_{j}q^{j}%
\]
using similar arguments we obtain%
\[
\gamma_{i-1}(q)=\overset{0}{\underset{j=-l}{\sum}}b_{j}q^{j}+\overset
{l}{\underset{j=1}{\sum}}b_{-j}q^{j},
\]
and so on. So we have computed the expansion of $G(T^{(i+1)})$ on the basis
$\{v_{\tau}\}$ and this expansion verifies $d_{\tau,T^{(i+1)}}(q)=0$ if
$\tau\vartriangleright T^{(i+1)}$. Finally notice that $\gamma_{s}%
(q)\in\mathbb{Z}[q,q^{-1}]$ by Proposition \ref{prop_A(T)}.

\begin{theorem}
Let $T\in\mathbf{OT}(n,\lambda).\;$Then $G(T)=\sum d_{\tau,T}(q)v_{\tau}$
where the coefficients $d_{\tau,T}(q)$ verify:

\textrm{(i)}: $d_{\tau,T}(q)\in\mathbb{Z}[q],$

\textrm{(ii)}: $d_{T,T}(q)=1$ and $d_{\tau,T}(0)=0$ for $\tau\neq T,$

\textrm{(iii)}: $d_{\tau,T}(q)\neq0$ only if $\tau$ and $T$ have the same
weight, and $\tau\trianglelefteq T$.
\end{theorem}

\noindent\textbf{Remark}: If we use $\mathcal{W}(\omega_{p})$ instead of
$W(\omega_{p}),$ the Lemma \ref{G_on_tabloids} is false.\ So we can not use a
similar triangular algorithm to compute $\{G(T)\}$ from $\{A(T)\}.$

\section{Example}

All the vectors occurring in our calculations are weight vectors.\ So we can
use our algorithm to compute the canonical basis of a single weight space. We
give below the matrix obtained for the $10$-dimension weight space of the
$U_{q}(so_{7})$-module $V(3,2,1)$ (i.e. $\lambda=\Lambda_{1}+\Lambda
_{2}+2\Lambda_{3}$) corresponding to the weight $(0,2,-1).$ Its columns and
rows are respectively labelled by the orthogonal tableaux $T_{i},$
$i=1,...,10$ and by the tabloids of weight $(0,2,-1)$ ordered from left to
right and top to bottom in increasing order for $\trianglelefteq$.\ The
orthogonal tableaux (which are tabloids) appear in the labelling of the rows
of the matrix.\ We have written them with bold font.

\noindent In this example we have, $G(T_{i})=A(T_{i})$ for $i\in
\{1,2,5,6,7,9,10\}$, $G(T_{3})=A(T_{3})-A(T_{2}),$ $G(T_{4})=A(T_{4})-\left(
q+\dfrac{1}{q}\right)  A(T_{2})$ and $G(T_{8})=A(T_{8})-A(T_{5}).$

\bigskip

{\tiny \hskip-16mm }%
\[
{\tiny
%\left[
%
%
%
%
%
%
%
%
%
%
%
%
%
%
%
%
%
%
%
%
%
%
%
%
%
%
%
%
%
%
%
%
%
%
%
%
%
%
%
%
%
%
%
%
%
%
%
%
%
%
\begin{array}
[c]{cccccccccccc}%
. &
\begin{array}
[c]{c}%
{222}\text{ }\\
{0\bar{3}{\ \ }}\\
{\bar{2}{\ }{\ }{\ }}%
\end{array}
{\ }\vspace{0.25cm} &
\begin{array}
[c]{c}%
{122}\text{ }\\
{0\bar{1}{\ \ }}\\
{\bar{3}{\ }{\ }{\ }}%
\end{array}
&
\begin{array}
[c]{c}%
{223}\text{ }\\
{0\bar{3}{\ \ }}\\
{\bar{3}{\ }{\ }{\ }}%
\end{array}
{\ } &
\begin{array}
[c]{c}%
{220}\text{ }\\
{3\bar{3}{\ \ }}\\
{\bar{3}{\ }{\ }{\ }}%
\end{array}
{\ } &
\begin{array}
[c]{c}%
{220}\text{ }\\
{0\bar{3}{\ \ }}\\
{\bar{3}{\ }{\ }{\ }}%
\end{array}
{\ } &
\begin{array}
[c]{c}%
{120}\text{ }\\
{2\bar{1}{\ \ }}\\
{\bar{3}{\ }{\ }{\ }}%
\end{array}
{\ } &
\begin{array}
[c]{c}%
{22\bar{3}}\text{ }\\
{30{\ \ }}\\
{\bar{3}{\ }{\ }{\ }}%
\end{array}
{\ } &
\begin{array}
[c]{c}%
{22\bar{3}}\text{ }\\
{3\bar{3}{\ \ }}\\
{0{\ }{\ }{\ }}%
\end{array}
{\ } &
\begin{array}
[c]{c}%
{12\bar{3}}\text{ }\\
{\bar{2}\bar{1}{\ \ }}\\
{0{\ }{\ }{\ }}%
\end{array}
{\ } &
\begin{array}
[c]{c}%
{12\bar{1}}\text{ }\\
{20{\ \ }}\\
{\bar{3}{\ }{\ }{\ }}%
\end{array}
{\ } &
\begin{array}
[c]{c}%
{12\bar{1}}\text{ }\\
{2\bar{3}{\ \ }}\\
{0{\ }{\ \ }}%
\end{array}
{\ }\\
{\
\begin{array}
[c]{c}%
221\text{ }\\
\bar{3}0\ \ \\
\bar{1}\ \ \
\end{array}
}\vspace{0.25cm} & q^{8} & . & . & . & . & . & . & . & . & . & .\\
{\
\begin{array}
[c]{c}%
221\text{ }\\
0\bar{3}\ \ \\
\bar{1}\ \ \
\end{array}
}\vspace{0.25cm} & q^{6} & q^{8} & . & . & . & . & . & . & . & . & .\\
{\
\begin{array}
[c]{c}%
221\text{ }\\
0\bar{1}\ \ \\
\bar{3}\ \ \
\end{array}
}\vspace{0.25cm} & . & q^{6} & . & . & . & . & . & . & . & . & .\\
{\
\begin{array}
[c]{c}%
012\text{ }\\
\bar{3}2\ \ \\
\bar{1}\ \ \
\end{array}
}\vspace{0.25cm} & . & . & . & q^{8} & . & . & . & . & . & . & .\\
{\
\begin{array}
[c]{c}%
212\text{ }\\
\bar{3}0\ \ \\
\bar{1}\ \ \
\end{array}
}\vspace{0.25cm} & q^{6} & . & . & q^{6} & . & . & q^{8} & . & . & . & .\\
{\
\begin{array}
[c]{c}%
212\text{ }\\
0\bar{3}\ \ \\
\bar{1}\ \ \
\end{array}
}\vspace{0.25cm} & q^{4} & q^{6} & . & . & . & . & q^{6} & . & . & q^{8} & .\\
{\
\begin{array}
[c]{c}%
212\text{ }\\
0\bar{1}\ \ \\
\bar{3}\ \ \
\end{array}
}\vspace{0.25cm} & . & q^{4} & . & . & . & . & . & . & . & q^{6} & .\\
{\
\begin{array}
[c]{c}%
122\text{ }\\
\bar{3}0\ \ \\
\bar{1}\ \ \
\end{array}
}\vspace{0.25cm} & q^{4} & . & . & . & . & . & q^{6} & . & q^{8} & . & .\\
{\
\begin{array}
[c]{c}%
222\text{ }\\
\bar{3}0\ \ \\
\bar{2}\ \ \
\end{array}
}\vspace{0.25cm} & q^{2} & . & . & q^{8}+q^{6} & q^{9}+q^{7} & . & q^{4} & . &
q^{6} & . & .\\
{\
\begin{array}
[c]{c}%
022\text{ }\\
00\ \ \\
\bar{3}\ \ \
\end{array}
}\vspace{0.25cm} & . & . & . & q^{5} & {q}^{6} & . & . & . & . & . & .\\
{\
\begin{array}
[c]{c}%
122\text{ }\\
0\bar{3}\ \ \\
\bar{1}\ \ \
\end{array}
}\vspace{0.25cm} & q^{2} & q^{4} & . & q^{6} & . & q^{8} & q^{4} & . & q^{6} &
q^{6} & .\\
{\
\begin{array}
[c]{c}%
\mathbf{222}\text{ }\\
\mathbf{0\bar{3}\ \ }\\
\mathbf{\bar{2}\ \ \ }%
\end{array}
}\vspace{0.25cm} & {1} & q^{2} & q^{8} & {q}^{6}+q^{4} & {q}^{7}+q^{5} & q^{6}%
& q^{2} & . & q^{4} & q^{4} & .\\
{\
\begin{array}
[c]{c}%
322\text{ }\\
0\bar{3}\ \ \\
\bar{3}\ \ \
\end{array}
}\vspace{0.25cm} & {.} & . & q^{6} & q^{4} & q^{5}-q^{9} & {.} & . & . & . &
. & .\\
{\
\begin{array}
[c]{c}%
022\text{ }\\
0\bar{3}\ \ \\
0\ \ \
\end{array}
}\vspace{0.25cm} & . & . & . & . & q^{4} & . & . & . & . & . & .\\
{\
\begin{array}
[c]{c}%
222\text{ }\\
0\bar{2}\ \ \\
\bar{3}\ \ \
\end{array}
}\vspace{0.25cm} & . & q^{2} & q^{4} & . & . & q^{6} & . & . & . & q^{4} & .\\
{\
\begin{array}
[c]{c}%
\mathbf{122}\text{ }\\
\mathbf{0\bar{1}\ \ }\\
\mathbf{\bar{3}\ \ \ }%
\end{array}
}\vspace{0.25cm} & . & 1 & . & . & {.} & q^{4} & . & . & . & q^{2} & .\\
{\
\begin{array}
[c]{c}%
232\text{ }\\
0\bar{3}\ \ \\
\bar{3}\ \ \
\end{array}
}\vspace{0.25cm} & . & . & q^{2} & q^{6}+q^{4} & q^{7}+q^{5} & q^{4} & . & . &
q^{8}+q^{6} & . & .\\
{\ }%
\begin{array}
[c]{c}%
202\text{ }\\
00\ \ \\
\bar{3}\ \ \
\end{array}
\vspace{0.25cm} & {.} & . & . & q^{3} & {q}^{4} & . & . & . & q^{5} & . & .\\
{\ }%
\begin{array}
[c]{c}%
102\text{ }\\
2\bar{3}\ \ \\
\bar{1}\ \ \
\end{array}
\vspace{0.25cm} & . & . & . & q^{4} & {.} & {q}^{6} & . & . & . & . & .\\
{\ }%
\begin{array}
[c]{c}%
202\text{ }\\
3\bar{3}\ \ \\
\bar{3}\ \ \
\end{array}
\vspace{0.25cm} & {.} & . & . & {q}^{2} & {q}^{5}+q^{3} & q^{4} & . & . &
q^{6}+q^{4} & . & .\\
{\ }%
\begin{array}
[c]{c}%
202\text{ }\\
0\bar{3}\ \ \\
0\ \ \
\end{array}
\vspace{0.25cm} & {.} & . & {.} & {.} & q^{2} & {.} & . & . & q^{3} & . & .
\end{array}
}%
\]
{\tiny $
%\right]
%
%
%
%
%
%
%
%
%
%
%
%
%
%
%
%
%
%
%
%
%
%
%
%
%
%
%
%
%
%
%
%
%
%
%
%
%
%
%
%
%
%
%
%
%
%
%
%
%
%
%
%
%
%
%
%
%
%
%
%
%
%
%
%
%
%
%
%
%
%
%
%
%
%
%
%
%
%
%
%
%
%
%
%
%
%
%
%
%
%
%
%
%
%
%
%
%
%
%
%
$ }

\bigskip

{\tiny \hskip-16mm $
%\left[
%
%
%
%
%
%
%
%
%
%
%
%
%
%
%
%
%
%
%
%
%
%
%
%
%
%
%
%
%
%
%
%
%
%
%
%
%
%
%
%
%
%
%
%
%
%
%
%
%
%
$}%
\[
{\tiny
\begin{array}
[c]{cccccccccccc}%
. &
\begin{array}
[c]{c}%
{222}\text{ }\\
{0\bar{3}{\ \ }}\\
{\bar{2}{\ }{\ }{\ }}%
\end{array}
{\ }\vspace{0.25cm} &
\begin{array}
[c]{c}%
{122}\text{ }\\
{0\bar{1}{\ \ }}\\
{\bar{3}{\ }{\ }{\ }}%
\end{array}
{\ } &
\begin{array}
[c]{c}%
{223}\text{ }\\
{0\bar{3}{\ \ }}\\
{\bar{3}{\ }{\ }{\ }}%
\end{array}
{\ } &
\begin{array}
[c]{c}%
{220}\text{ }\\
{3\bar{3}{\ \ }}\\
{\bar{3}{\ }{\ }{\ }}%
\end{array}
{\ } &
\begin{array}
[c]{c}%
{220}\text{ }\\
{0\bar{3}{\ \ }}\\
{\bar{3}{\ }{\ }{\ }}%
\end{array}
{\ } &
\begin{array}
[c]{c}%
{120}\text{ }\\
{2\bar{1}{\ \ }}\\
{\bar{3}{\ }{\ }{\ }}%
\end{array}
{\ } &
\begin{array}
[c]{c}%
{22\bar{3}}\text{ }\\
{30{\ \ }}\\
{\bar{3}{\ }{\ }{\ }}%
\end{array}
{\ } &
\begin{array}
[c]{c}%
{22\bar{3}}\text{ }\\
{3\bar{3}{\ \ }}\\
{0{\ }{\ }{\ }}%
\end{array}
{\ } &
\begin{array}
[c]{c}%
{12\bar{3}}\text{ }\\
{\bar{2}\bar{1}{\ \ }}\\
{0{\ }{\ }{\ }}%
\end{array}
{\ } &
\begin{array}
[c]{c}%
{12\bar{1}}\text{ }\\
{20{\ \ }}\\
{\bar{3}{\ }{\ }{\ }}%
\end{array}
{\ } &
\begin{array}
[c]{c}%
{12\bar{1}}\text{ }\\
{2\bar{3}{\ \ }}\\
{0{\ }{\ \ }}%
\end{array}
{\ }\\
{\
\begin{array}
[c]{c}%
102\text{ }\\
2\bar{1}\ \ \\
\bar{3}\ \ \
\end{array}
}\vspace{0.25cm} & . & . & . & . & . & q^{2} & . & . & q^{4} & . & .\\
{\
\begin{array}
[c]{c}%
1\bar{3}2\text{ }\\
2\bar{1}\ \ \\
0\ \ \
\end{array}
}\vspace{0.25cm} & . & . & . & . & . & . & . & . & q^{2} & . & .\\
{\
\begin{array}
[c]{c}%
\mathbf{223}\text{ }\\
\mathbf{0\bar{3}\ \ }\\
\mathbf{\bar{3}\ \ \ }%
\end{array}
}\vspace{0.25cm} & . & . & 1 & q^{4}+q^{2} & q^{5}+q^{3} & q^{2} & q^{6} &
q^{8} & q^{6}+q^{4} & q^{4} & .\\
{\ }%
\begin{array}
[c]{c}%
210\text{ }\\
\bar{3}2\ \ \\
\bar{1}\ \ \
\end{array}
\vspace{0.25cm} & . & . & . & q^{4} & . & . & q^{6} & . & . & . & .\\
{\ }%
\begin{array}
[c]{c}%
220\text{ }\\
00\ \ \\
\bar{3}\ \ \
\end{array}
\vspace{0.25cm} & . & . & . & . & q^{2} & . & q^{3} & q^{5} & q^{3} & . & .\\
{\ }%
\begin{array}
[c]{c}%
120\text{ }\\
2\bar{3}\ \ \\
\bar{1}\ \ \
\end{array}
\vspace{0.25cm} & . & . & . & q^{2} & . & q^{4} & q^{4} & . & q^{6} & q^{6} &
q^{8}\\
{\ }%
\begin{array}
[c]{c}%
\mathbf{220}\text{ }\\
\mathbf{3\bar{3}\ \ }\\
\mathbf{\bar{3}\ \ \ }%
\end{array}
\vspace{0.25cm} & . & . & . & 1 & q^{3}+q & q^{2} & q^{2} & q^{6}+q^{4} &
2q^{4}+q^{2} & q^{4} & q^{6}\\
{\ }%
\begin{array}
[c]{c}%
\mathbf{220}\text{ }\\
\mathbf{0\bar{3}\ \ }\\
\mathbf{0\ \ \ }%
\end{array}
\vspace{0.25cm} & . & . & . & . & 1 & . & . & q^{3} & q & . & .\\
{\ }%
\begin{array}
[c]{c}%
\mathbf{120}\text{ }\\
\mathbf{2\bar{1}\ \ }\\
\mathbf{\bar{3}\ \ \ }%
\end{array}
\vspace{0.25cm} & . & . & . & . & {.} & 1 & . & . & q^{2} & q^{2} & q^{4}\\
{\ }%
\begin{array}
[c]{c}%
21\bar{3}\text{ }\\
02\ \ \\
\bar{1}\ \ \
\end{array}
\vspace{0.25cm} & . & . & . & . & . & . & q^{4} & . & . & q^{6} & .\\
{\ }%
\begin{array}
[c]{c}%
22\bar{3}\text{ }\\
03\ \ \\
\bar{3}\ \ \
\end{array}
\vspace{0.25cm} & {.} & . & . & {.} & {.} & . & q^{2} & q^{4} & . & q^{4} &
.\\
{\ }%
\begin{array}
[c]{c}%
12\bar{3}\text{ }\\
20\ \ \\
\bar{1}\ \ \
\end{array}
\vspace{0.25cm} & {.} & . & . & . & . & {.} & q^{2} & . & q^{4} & q^{4} &
q^{6}\\
{\ }%
\begin{array}
[c]{c}%
\mathbf{22\bar{3}}\text{ }\\
\mathbf{30\ \ }\\
\mathbf{\bar{3}\ \ \ }%
\end{array}
\vspace{0.25cm} & . & . & . & . & . & . & 1 & q^{4}+q^{2} & q^{2} & q^{2} &
q^{4}\\
{\ }%
\begin{array}
[c]{c}%
22\bar{3}\text{ }\\
00\ \ \\
0\ \ \
\end{array}
\vspace{0.25cm} & . & . & . & . & . & . & . & q & . & . & .\\
{\ }%
\begin{array}
[c]{c}%
\mathbf{22\bar{3}}\text{ }\\
\mathbf{3\bar{3}\ \ }\\
\mathbf{0\ \ \ }%
\end{array}
\vspace{0.25cm} & . & . & . & . & {.} & . & . & 1 & q^{2} & . & q^{4}\\
{\ }%
\begin{array}
[c]{c}%
\mathbf{12\bar{3}}\text{ }\\
\mathbf{2\bar{1}\ \ }\\
\mathbf{0\ \ \ }%
\end{array}
\vspace{0.25cm} & . & . & . & . & . & . & . & . & 1 & . & q^{2}\\
{\ }%
\begin{array}
[c]{c}%
21\bar{1}\text{ }\\
02\ \ \\
\bar{3}\ \ \
\end{array}
\vspace{0.25cm} & {.} & . & . & . & {.} & . & . & . & . & q^{2} & .\\
{\ }%
\begin{array}
[c]{c}%
\mathbf{12\bar{1}}\text{ }\\
\mathbf{20\ \ }\\
\mathbf{\bar{3}\ \ \ }%
\end{array}
\vspace{0.25cm} & . & . & . & . & {.} & {.} & . & . & . & 1 & q^{2}\\
{\ }%
\begin{array}
[c]{c}%
\mathbf{12\bar{1}}\text{ }\\
\mathbf{2\bar{3}\ \ }\\
\mathbf{0\ \ \ }%
\end{array}
\vspace{0.25cm} & {.} & . & . & {.} & {.} & . & . & . & . & . & 1
\end{array}
}%
\]
{\tiny $
%\right]
%
%
%
%
%
%
%
%
%
%
%
%
%
%
%
%
%
%
%
%
%
%
%
%
%
%
%
%
%
%
%
%
%
%
%
%
%
%
%
%
%
%
%
%
%
%
%
%
%
%
%
%
%
%
%
%
%
%
%
%
%
%
%
%
%
%
%
%
%
%
%
%
%
%
%
%
%
%
%
%
%
%
%
%
%
%
%
%
%
%
%
%
%
%
%
%
%
%
%
%
$ }

\end{document}